\documentclass[11pt]{sat}

\newif\ifsattoc\sattoctrue
\newread\testfl\immediate\openin\testfl=\jobname.tic
    \ifeof\testfl\sattocfalse\fi\closein\testfl

\newcommand{\sect}[1]{\section{#1}\setcounter{equation}{0}}

\title{Logarithmic Potential Theory with \\ Applications to
Approximation Theory}
\def\shorttitle{Logarithmic Potential Theory}

\author{E.B. Saff \footnote{Research was
supported, in part, by the U. S. National Science Foundation under grant
DMS-0808093.}}

\def\shortauthor{E. B. Saff}

\def\versiondate{12 October 2010}

\def\abstracttext{
We provide an introduction to logarithmic potential theory in the complex plane
that particularly emphasizes its usefulness in the theory of polynomial and
rational approximation.  The reader is invited to  explore the notions of
Fekete points, logarithmic capacity, and Chebyshev constant through a variety
of examples and exercises. Many of the fundamental theorems of potential
theory, such as Frostman's theorem, the Riesz Decomposition Theorem,
the Principle of Domination, etc., are given along with essential ideas for
their proofs. Equilibrium measures and potentials and their connections with
Green functions and conformal mappings are presented.  Moreover, we discuss
extensions of the classical potential theoretic results to the case when an
external field is present.
}

\def\MSCnumbers{ Primary: 31Cxx, 41Axx} 

\def\keywords{Logarithmic potential, Polynomial approximation, Rational
approximation, Transfinite diameter, Capacity, Chebyshev constant, Fekete
points, Equilibrium potential, Superharmonic functions, Subharmonic functions,
Green functions, Rates of polynomial and rational approximation, Condenser
capacity,  External fields.}

\def\dword#1{{\bf #1}}
\def\dd{\,{\rm d}}  
\def\ee{{\rm e}}  
\def\ii{{\rm i}}  

\def\startpagenumber{165}
\def\volumenumber{5}
\def\year{2010}

\newcommand{\beginddoc}{
\begin{document}
\maketitle
\begin{abstract}
\abstracttext
\vskip1pt MSC: \MSCnumbers
\ifx\keywords\empty\else\vskip1pt Keywords: \keywords\fi
\end{abstract}
\insert\footins{\scriptsize
\medskip
\baselineskip 8pt
\leftline{Surveys in Approximation Theory}
\leftline{Volume \volumenumber, \year.
pp.~\thepage--\pageref{endpage}.}
\leftline{\copyright\ \year\ Surveys in Approximation Theory.}
\leftline{ISSN 1555-578X}
\leftline{All rights of reproduction in any form reserved.}
\smallskip
\par\allowbreak}
\ifsattoc\else\tableofcontents\fi}
\renewcommand\rightmark{\ifodd\thepage{\it \hfill\shorttitle\hfill}\else {\it \hfill\shortauthor\hfill}\fi}
\markboth{{\it \shortauthor}}{{\it \shorttitle}}
\markright{{\it \shorttitle}}
\def\endddoc{\label{endpage}\end{document}}
\date{{\small \versiondate}}
\setlength\oddsidemargin{0pc}
\setlength\evensidemargin{0pc}
\setlength\topmargin{0in}
\setlength\textwidth{6.5in}
\setlength\textheight{8.6in}
\usepackage{amssymb} 
\usepackage{graphicx} \usepackage{amsmath} \usepackage{amsfonts}
\usepackage{enumerate} \usepackage{verbatim} \usepackage{footmisc}

\newcommand{\C}{\mathbb C}\newcommand{\R}{\mathbb R}
\newcommand{\M}{{\cal M}}\newcommand{\F}{{\cal F}}
\newcommand{\Po}{{\cal P}}\newcommand{\A}{{\cal A}}
\newcommand{\cp}{\mbox{\textnormal{cap}}}
\newcommand{\cheb}{\mbox{\textnormal{cheb}}}
\renewcommand{\thefootnote}{\fnsymbol{footnote}}
\beginddoc
\ifsattoc
\bigskip
\def\toczer{0}\def\tochalf{.5}\def\toctwo{1}
\def\tocindent{0}
\def\numberline#1{\hskip\tocindent truecm{} #1\hskip1em}
\newread\testfl
\def\inputifthere#1{\immediate\openin\testfl=#1
    \ifeof\testfl\message{(#1 does not yet exist)}
    \else\input#1\fi\closein\testfl}
\countdef\counter=255
\def\diamondleaders{\global\advance\counter by 1
  \ifodd\counter \kern-10pt \fi
  \leaders\hbox to 15pt{\ifodd\counter \kern13pt \else\kern3pt \fi
  \hss.\hss}\hfill}
\newdimen\lextent
\newtoks\writestuff
\medskip
\begingroup
\small
\def\contentsline#1#2#3{
\ifx#1{section}\let\tocindent\toczer\else\let\tocindent\tochalf\fi
\setbox1=\hbox{#2}\ifnum\wd1>\lextent\lextent\wd1\fi}
\lextent0pt\inputifthere{\jobname.tic}\advance\lextent by 2em\relax
\def\contentsline#1#2#3{
\ifx#1{section}\let\tocindent\toczer\else\let\tocindent\tochalf\fi
\writestuff={#2}
\centerline{\hbox to \lextent{\rm\the\writestuff%
\ifx\empty#3\else\diamondleaders{}
\hfil\hbox to 2 em\fi{\hss#3}}}}
\input \jobname.tic\endgroup \newpage 
\fi




%


\newcounter{example} \setcounter{example}{1}
\setcounter{section}{-1}

\sect{Introduction} Logarithmic potential theory is an elegant blend
of real and complex analysis that has had a profound effect on many
recent developments in approximation theory.  Since logarithmic
potentials have a direct connection with polynomial and rational
functions, the tools provided by classical potential theory and its
extensions to cases when an external field (or weight) is present, have
resolved some long-standing problems concerning orthogonal polynomials,
rates of polynomial and rational approximation, convergence behavior of
Pad\'{e} approximants (both classical and multi-point), to name but a
few.  Here are some problems where potential theory has played a crucial
role:

\begin{enumerate}[(i)]
\item {\em Rate of Polynomial Approximation}:
Let $f$ be analytic on a compact set $E$ of the complex plane $\C$,
whose complement $\overline{\C} \setminus E$ is connected.  How well can
$f$ be uniformly approximated on $E$ by polynomials (in $z$) of degree
at most $n$?
\item {\em Asymptotic Behavior of Zeros of Polynomials}:
Let $p^*_n(x)$ denote the polynomial of degree at most $n$ of best
uniform approximation to a continuous function on $[-1,1]$, say $f(x) =
|x|$.  In the complex plane, $p^*_n$ has $n$ zeros (at
most).\footnote{By symmetry and uniqueness of the best approximants,
$p^*_{2n+1}=p^*_{2n}$ for $f(x)=|x|$.}  Where are these zeros located as
$n \rightarrow \infty$?
\item {\em Fast Decreasing Polynomials}:  Given
$\varphi \in C[-1,1]$, does there exist a sequence of ``needle-like''
polynomials $(p_n)$, $\deg p_n \le n$, such that $p_n(0) = 1$ and
$|p_n(x)|\le C\ee^{-cn\varphi(x)}, \ x\in [-1,1]$, for some positive
constants $c,C$?
\item {\em Recurrence Coefficients for Orthogonal Polynomials}:
Let $\{p_n\}$ denote orthonormal polynomials with respect
to the weight $\exp\left(-|x|^{\alpha}\right), \ \alpha>0$, on $\R$.
That is, $$\int_{-\infty}^{\infty} p_m(x)p_n(x)\ee^{-|x|^{\alpha}}\dd x =
\delta_{mn}.$$ Then the $p_n$ satisfy a 3-term recurrence of the form
$$x p_n(x) = a_{n+1} p_{n+1}(x) + a_{n} p_{n-1}(x), \ \ \ \ \ n = 1, 2,
\ldots,$$ where $(a_n)$ is the sequence of recurrence coefficients.
These coefficients go to infinity as $n$ increases, but exactly what is
their asymptotic growth rate?
\item {\em Generalized Weierstrass Problem}:
A famous theorem of Weierstrass states that ${f \in C[-1,1]}$
if and only if there exists a sequence of polynomials $(p_n)$, $\deg
p_n \le n$, such that $p_n \rightarrow f$ uniformly on $[-1,1]$.  But
how would you characterize those ${f \in C[-1,1]}$ that are uniform
limits on $[-1, 1]$ of ``incomplete'' polynomials of the form $q_{2n}(x)
= \sum_{k=n}^{2n} a_k x^k$, for which half the coefficients are missing?
More generally, what functions $f$ are the uniform limits of {\em
weighted polynomials} of the form $w(x)^n p_n(x)$, where the power of
the weight matches the degree of the polynomial?
\item {\em Optimal Point Arrangements on the Sphere}:
How well separated are $N\ (\ge2)$
points on the unit sphere $S^2 = \{ x\in\R^3 : |x|=1 \}$ that maximize
the product of their pairwise distances: $$\prod_{1\le i<j \le N} |x_i -
x_j| \ \ {?}$$
\item {\em Rational Approximation}:  Determine the
rate of best uniform approximation to $\ee^{-x}$ on $[0, +\infty)$ by
rational functions of the form $p_n(x)/q_n(x), \ \deg p_n\le n, \ \deg
q_n \le n$.
\end{enumerate}

In this article we provide an introduction to the tools of classical
and ``weighted'' potential theory that are the keys to resolving the
above questions.  The essential reason for the usefulness of
potential theory in obtaining results on polynomials is the fact that
for any monic polynomial $p(z) = \prod_{k=1}^n (z-z_k)$ the function
$\log(1/|p(z)|)$ can be written as a logarithmic potential:
$$\log\frac{1}{|p(z)|} = \int \log\frac{1}{|z-t|}\dd\nu(t),$$
where $\nu$ is the discrete measure with mass $1$ at each of the zeros of $p$.

\sect{Transfinite Diameter, Capacity, and Chebyshev Constant}\label{section classical}

We begin by introducing three ``different'' quantities associated with a
compact (closed and bounded) set in the plane.
\newline

\noindent {\bf A Geometric Problem.}  Place $n$ points on a compact set
$E$ so that they are ``as far apart'' as possible in the sense of the
{\em geometric mean} of the pairwise distances between the points.
Since the number of different pairs of $n$ points is $n(n-1)/2$, we
consider the quantity
\begin{equation}
\delta_n(E):=\max_{z_1,\ldots,z_n\in E}\left(\prod_{1\leq i<j\leq
n}|z_i-z_j|\right)^{2/n(n-1)}.  \label{geom prob}
 \end{equation}
Any system of points $\F_n:=\left\{z_1^{(n)},\ldots,z_n^{(n)}\right\}$ for
which the maximum is attained is called an \dword{$n$-point Fekete set}
for $E$; the points $z_i^{(n)}$ in $\F_n$ are called \dword{Fekete
points}.

For example, if $n=2$, then $\F_2=\left\{z_1^{(2)},z_2^{(2)}\right\}$,
where $\left|z_1^{(2)}-z_2^{(2)}\right|=\mbox{ diam }E$. Obviously, any
such points lie on the boundary of $E$.  In general, it follows from the
maximum modulus principle for analytic functions that for all $n$, the
Fekete sets lie on the {\em outer boundary} $\partial_{\infty}E$, that
is, the boundary of the unbounded component of the complement of $E$.
\newline
\newline
{\bf \em Exercise.}   Prove that the determinant of
the $n\times n$ Vandermonde matrix $[z_i^j], \  1\le i \le n, \ 0 \le j
\le n-1$, is given by $\prod_{1\le i<j\le n} (z_j - z_i)$.
Consequently, an $n$-point Fekete set for $E$ maximizes the modulus of
this determinant over all $n$-point subsets of $E$.
\newline
\newline
{\bf \em Exercise.}  Let $E$ be the closed unit disk (or the unit
circle).  Prove that the set of $n$th roots of unity is an $n$-point
Fekete set for $E$ (and so is any of its rotations) and that
$\delta_n(E) = n^{1/(n-1)}$.  [Hint: Use Hadamard's inequality for
determinants.]
\newline

If $E=[-1, 1]$, then the set $\F_n$ turns out to be unique and it
coincides with the zeros of $(1-x^2)P'_{n-1}(x)$, where $P_{n-1}$ is the
Legendre polynomial of degree $n-1$ (cf.  [Sz]).  \par Fekete points are
``good points'' for polynomial interpolation.  We denote by
$\mathcal{P}_n$ the linear space of all algebraic polynomials with
complex coefficients of degree at most $n$.  Recall that if $z_1,
\ldots, z_{n+1}$ are any $n+1$ distinct points, then the unique
polynomial in $\mathcal{P}_n$ that interpolates a function $f$ in these
points is given by $$p_n(z) = \sum_{k=1}^{n+1} f(z_k) L_k(z),$$ where
$L_k(z)$ is the \dword{fundamental Lagrange polynomial} that satisfies
$L_k(z_j) = \delta_{jk}$.
\newline
\newline
{\bf \em Exercise.}  Prove
that if $\{z_1, \ldots, z_{n+1}\}$ is an $(n+1)$-point Fekete set for a
compact set $E$, then the associated fundamental Lagrange polynomials
satisfy $|L_k(z)| \le 1$ for all $z\in E$.  Furthermore, show that if
$P_n \in \mathcal{P}_n$, then $$\lVert P_n \rVert_E \le (n+1) \lVert P_n
\rVert_{\F_{n+1}},$$ where $\lVert \cdot \rVert_A$ denotes the $\sup$
norm on $A$.
\newline
\newline
{\bf \em Exercise.}   Prove that if
$P_n(f;z)$ denotes the polynomial of degree at most $n$ that
interpolates a continuous function $f$ in an $(n+1)$-point Fekete set
for $E$, then $$\lVert f - P_n(f;\cdot) \rVert_E \le (n+2) \lVert f -
p_n^* \rVert_E,$$ where $p_n^*$ is the best uniform approximation to $f$
on $E$ out of $\mathcal{P}_n$.
\newline

\par On taking the logarithm in (\ref{geom prob}), we see that the max
problem in (\ref{geom prob}) is equivalent to the {\em minimization
problem}
\begin{equation} \mathcal{E}_n(E) := \min_{z_1, \ldots, z_n \in
E} \sum_{1\le i<j\le n} \log\frac{1}{|z_i - z_j|}.  \label{min prob}
\end{equation}
The summation in (\ref{min prob}) can be interpreted as
the energy of a system of $n$ like-charged particles located at the
points $\{z_i\}_{i=1}^n$, where the repelling force between two
particles is proportional to the reciprocal of the distance between
them.  Thus
\begin{equation} \mathcal{E}_n(E) = \frac{n(n-1)}{2} \log
\frac{1}{\delta_n(E)} \label{log energy}
 \end{equation}
denotes the minimal logarithmic energy that can be attained by $n$ particles
that are constrained to lie on $E$.  Any set of $n$ points that attains this
minimal energy is called an \dword{equilibrium configuration} for $E$;
that is, a Fekete set $\F_n$ represents an $n$-point equilibrium
configuration for $E$.

Essential questions are:
\begin{enumerate}[(i)]
\item What is the
asymptotic behavior of the minimal energy $\mathcal{E}_n(E)$ (or,
equivalently, of $\delta_n(E)$) as $n\rightarrow \infty$?
\item How are
optimal configurations (Fekete points) distributed on $E$ as
$n\rightarrow \infty$?
\end{enumerate}
As a first step we establish
\newline
\newline
{\bf Lemma \arabic{section}.\arabic{example}.}\addtocounter{example}{1}
{\sl The sequence $\left( \frac{\mathcal{E}_n(E)}{n(n-1)}
\right)_{n=2}^{\infty}$ is increasing (i.e., nondecreasing) or,
equivalently, the sequence $( \delta_n(E) )_{n=2}^\infty$ is
decreasing (i.e., nonincreasing).}
\newline
\newline
{\bf Proof.}  With
$\F_n = \left\{z_k^{(n)}\right\}_{k=1}^n$ we have, for each
$k=1,\ldots,n$,
\begin{equation*} \begin{split} \mathcal{E}_n(E) & =
\sum_{i \ne k} \log \frac{1}{\left|z_i^{(n)} - z_k^{(n)}\right|} +
\sum_{\substack{1\le i<j \le n \\ i\ne k \\ j\ne k}} \log
\frac{1}{\left|z_i^{(n)} - z_j^{(n)}\right|} \\ & \ge \sum_{i \ne k}
\log \frac{1}{\left|z_i^{(n)} - z_k^{(n)}\right|} +
\mathcal{E}_{n-1}(E).
  \end{split}
 \end{equation*}
Now add these $n$ inequalities together and divide by $n(n-1)(n-2)$ to get
result.  $\Box$
\newline
\newline

The sequence $( \delta_n(E) )$ therefore has a limit\footnote{If $E$
consists of only finitely many points, then $\tau(E) = 0$.  Why?}
\begin{equation} \boxed{\tau(E) := \lim_{n\rightarrow \infty}
\delta_n(E),} \label{transfinite diameter}
 \end{equation} which is
called the \dword{transfinite diameter} of $E$.  For example, the
transfinite diameter of the disk $E= \{ z \in \C : |z| \le R \}$ is $R$
since $\delta_n(E) = Rn^{1/(n-1)} \rightarrow R$ as
$n\rightarrow\infty$.

Note that $0 \le \tau(E) \le \mbox{ diam } E$ and that $E_1 \subset E_2$
implies $\tau(E_1) \le \tau(E_2)$.
\newline
\newline
{\bf \em Exercise.}
Let $aE +b  := \{ az+b : z\in E \}$, with $a,b$ fixed complex constants.
Prove that $\tau(aE+b) = |a| \tau(E)$ for any compact set $E\subset \C$.
\newline
\newline
 {\bf \em Exercise.}   Show that the closed set $E =
\{0\} \cup \{1/k : k = 1, 2, \ldots \}$ has transfinite diameter zero.
\newline
\newline
{\em Remark.}  The transfinite diameter $\tau$
(considered as a set function) has some of the properties of Lebesgue
measure on compact subsets of $\C$;  in fact, if $E$ is the closed
interval $[a,b]$, then $\tau([a,b]) = (b-a)/4$.  However, $\tau$ {\em
fails to be subadditive}; $\tau(E_1\cup E_2)$ may exceed the sum
$\tau(E_1) + \tau(E_2)$.

To investigate the asymptotic behavior of a sequence of Fekete sets
$\F_n, ~ n=2,3,\ldots$, we utilize {\em weak-star convergence} of
measures.

{\bf Definition \arabic{section}.\arabic{example}.}\addtocounter{example}{1}
Let $\mu_n$ be a sequence of finite positive measures with
supports\footnote{Recall
that a point $z_0$ belongs to $\mbox{supp}(\mu)$ if and only if every
open set containing $z_0$ has positive $\mu$-measure.}
$\mbox{supp}(\mu_n) \subset K$ for all $n$, where $K$ is some compact
set.  We write $\mu_n \stackrel{*}{\rightarrow} \mu$ if
\begin{equation}
\lim_{n\rightarrow\infty} \int f \dd\mu_n = \int  f \dd\mu \ \ \ \ \forall
f \in C(K).  \label{weak}
 \end{equation}
(If $\mu_n (K) \le M$ for some
constant $M$ and all $n$ (which clearly holds when $\mu$ is a finite
measure), this is equivalent to pointwise convergence in the dual space
of $C(K)$.)  The same definition applies to signed measures and complex
measures.  In (\ref{weak}), we can always take $K$ to be the extended
complex plane $\overline{\C}$; however, knowing a specific compact set
$K$ that contains all the supports of the $\mu_n$'s serves to remind us
that the limit measure will also be supported on $K$.
\newline

For a discrete set consisting of $n$ points of $\C$, say $A_n = \{z_1,
\ldots, z_n\}$, we associate the \dword{normalized counting measure}
$$\nu (A_n) := \frac{1}{n}\sum_{k=1}^n \delta_{z_k},$$
where $\delta_z$ is the unit point mass at $z$.
\newline
\newline
{\bf Example \arabic{section}.\arabic{example}.} \addtocounter{example}{1}  If
$A_{n+1}$ consists of the $n+1$ \dword{Chebyshev nodes} for $[-1, 1]$;
that is $A_{n+1} = \{ \cos(k\pi/n) : k = 0, 1, \ldots, n\}$, then
\begin{equation} \nu (A_n) \stackrel{*}{\rightarrow} \frac{\dd x}{\pi
\sqrt{1-x^2}}, \label{arcsin}
 \end{equation} which is the \dword{arcsine
distribution} on $[-1,1]$.  The nodes $A_{n+1}$ are the extreme points
of the \dword{Chebyshev polynomials} $T_n(x) = \cos(n \arccos x)$, which
are orthogonal on $[-1, 1]$ with respect to the arcsine distribution.
Verify (\ref{arcsin})!
\newline
\newline
{\bf Example \arabic{section}.\arabic{example}.} \addtocounter{example}{1}
If $A_n$ consists of the $n$th roots of unity, then
$\nu(A_n)\stackrel{*}{\rightarrow} \frac{1}{2\pi} \dd\theta$, where
$\dd\theta$ is arclength on the unit circle $|z|=1$.
\newline
\newline
{\bf \em Exercise.} Let $\lambda$ be a real {\em irrational} number and
let $A_n := \{ \exp (\lambda k\pi\ii) : k=1,\ldots, n\}$.  Prove that
$\nu(A_n)\stackrel{*}{\rightarrow} \frac{1}{2\pi} \dd\theta$.  What
happens if $\lambda$ is rational?
\newline

As we shall see, many of the results of potential theory are formulated
for semi-continuous functions.
\newline
\newline
{\bf Definition \arabic{section}.\arabic{example}.} \addtocounter{example}{1}
A function $f : D \rightarrow (-\infty, \infty]$ ($f$ omits the value
$-\infty$) is \dword{lower semi-continuous} (\dword{l.s.c.}) on the set $D\subset \C$ if
it satisfies any of the following equivalent conditions:
\begin{enumerate}[(i)]
\item $\{ z \in D : f(z)>\alpha \}$ is open
relative to $D$ for every $\alpha \in \R$;
\item For every $z_0 \in D$,
$$f(z_0) \le \liminf_{z\rightarrow z_0} f(z);$$
\item For every compact
subset $K\subset D$, there exists an increasing sequence of continuous
functions on $K$ with pointwise limit $f$.
  \end{enumerate}

\noindent {\bf \em Exercise.}  Prove that if $f$ is l.s.c.\  on a compact
set $K$, then $f$ attains its minimum on $K$.
\newline

Important for us is the fact, which follows from property (iii) and the
Monotone Convergence Theorem, that if $f$ is l.s.c.\  on a compact set
$K$, then
\begin{equation} \int_K f \dd\mu \le
\liminf_{n\rightarrow\infty} \int_K f \dd\mu_n \label{monotone
convergence}
 \end{equation}
wherever $\mu_n\stackrel{*}{\rightarrow}\mu$
and $\mbox{supp}(\mu_n) \subset K$ for all $n$.
\newline
\newline
{\bf \em Exercise.}  Prove that if $\mu_n \stackrel{*}{\rightarrow} \mu$,
then for any bounded Borel set $E$, $$\mu(\mathring{E}) \le
\liminf_{n\rightarrow\infty} \mu_n(E) \le \limsup_{n\rightarrow\infty}
\mu_n(E) \le \mu(\overline{E}),$$ where $\mathring{E}$ and
$\overline{E}$ denote, respectively, the interior of $E$ and the closure
of $E$.  [Hint: First show that the characteristic function of an open
set is l.s.c.]
\newline

Our goal now is to determine the weak-star limit (if it exists) for the
sequence of normalized counting measures $\nu(\F_n)$ in the Fekete
points for a given compact set $E$.  For this purpose we study the
continuous analogue of the discrete minimum energy problem (\ref{min
prob}).
\newline
\newline
{\bf Electrostatics Problem for a Capacitor.}
Place a unit positive charge on a compact set $E$ so that equilibrium is
attained in the sense that energy is minimized.  Again it is assumed
that the repulsive force between like-charged particles located at
points $z$ and $t$ is proportional to $1/|z-t|$.

To create a mathematical framework for this problem, we let $\M(E)$
denote the collection of all positive unit Borel measures $\mu$
supported on $E$ (so that $\M(E)$ contains all possible distributions of
charges placed on $E$).  The \dword{logarithmic potential} associated with
$\mu$ is
\begin{equation} U^\mu(z):=\int\log\frac{1}{|z-t|}\dd\mu(t),
\label{log potential}
\end{equation}
 which is harmonic outside the
support $\mbox{supp}(\mu)$  of $\mu$ and is l.s.c.\  in $\C$ since
$$U^\mu(z) = \lim_{M\rightarrow\infty} \int \min\left(M,\ \log
\frac{1}{|z-t|}\right) \dd\mu(t).$$ The \dword{energy} of such a potential
is defined by
\begin{equation} I(\mu):=\int U^\mu
\dd\mu=\int\int\log\frac{1}{|z-t|}\dd\mu(t)\dd\mu(z).  \label{energy of
potential}
\end{equation}
Thus, the electrostatics problem involves the
determination of \begin{equation} V_E:=\inf\{I(\mu): \; \mu\in\M(E)\},
\label{electrostatics problem}
 \end{equation}
which is called the \dword{Robin constant} for $E$.
Note that since $E$ is bounded, we have
$$-\infty<V_E\leq+\infty.$$

First we establish the existence of a measure $\mu_E \in \M(E)$ for
which the ``inf'' is attained.  For this purpose we use
\newline
\newline
 {\bf Lemma
\arabic{section}.\arabic{example}}\addtocounter{example}{1} (\dword{Principle
of Descent}).  {\sl Let $\mu_n$ be a sequence of measures in
$\M(E)$ that converges weak-star to some $\mu \in \M(E)$.
Then, for all $z \in \C$,
\begin{equation} U^{\mu}(z) \le
\liminf_{n\rightarrow\infty} U^{\mu_n}(z), \label{descent U}
\end{equation}
and, furthermore, \begin{equation} I(\mu) \le
\liminf_{n\rightarrow\infty} I(\mu_n).  \label{descent I}
 \end{equation}
 }
\newline
 {\bf Proof.}  Inequality (\ref{descent U}) follows from
(\ref{monotone convergence}) on observing that $\log 1/|z-t|$ is l.s.c.
in $t$.  Inequality (\ref{descent I}) follows similarly, on observing
that $\mu_n \times \mu_n$ converges weak-star to $\mu \times \mu$.
$\Box$
\newline
\newline
{\bf Lemma \arabic{section}.\arabic{example}.}\addtocounter{example}{1}
{\sl There is some $\mu_E \in \M(E)$ such that $I(\mu_E) = V_E$.}
\newline
\newline
{\bf Proof.}  By the Banach-Alaoglu Theorem, $\M(E)$ is compact in the
weak-star topology (this fact is also known as \dword{Helly's Selection
Theorem}).  Let $\mu_n$ be a sequence in $\M(E)$ satisfying
$\lim_{n\rightarrow\infty} I(\mu_n) = V_E$ and let $\nu$ denote some
weak-star cluster point of the $\mu_n$.  Then, by the Principle of
Descent and the definition of $V_E$, we obtain $V_E = I(\nu)$.  $\Box$
\newline

When $V_E = +\infty$ (for example this is the case if $E$ is countable),
then every measure $\mu\in\M(E)$ is a minimizing measure.  However, if
$V_E$ is finite, it follows from the strict convexity of $I(\mu)$ on
$\M(E)$ (cf.~[ST]) that there exists a {\em unique} measure $\mu_E$
such that $V_E = I(\mu_E)$.  In this case, we call $\mu_E$ the
\dword{equilibrium measure} for $E$, and $U^{\mu_E}$ the \dword{equilibrium {\rm
or} conductor potential} for $E$.
\newline
\newline
{\bf Definition \arabic{section}.\arabic{example}.}\addtocounter{example}{1}
The \dword{logarithmic capacity} of $E$, denoted by $\cp(E)$, is defined by
\begin{equation} \boxed{\cp(E) := \ee^{-V_E}.} \label{log cap}
\end{equation}
If $V_E=+\infty$, we set $\cp(E) = 0$; such sets $E$ are
called \dword{polar sets} because they correspond to sets where potentials
can equal $+\infty$.  More generally, an arbitrary set $E \subset
\mathbb{C}$ is said to be \dword{polar} if every closed subset of $E$ is polar.
In electrostatic terms, polar sets are ``too small'' to hold a
charge.\footnote{Somewhat surprising is the fact that the classical
``$1/3$ Cantor set,'' which has 1-dimensional Lebesgue measure zero, has
{\em positive} capacity.  The precise value of this capacity is as yet
still unknown (see
 [R2]  for some numerical approximation methods).}
\newline
\newline
{\bf \em Exercise.}  Prove that any countable set $E$ is a polar set.
\newline

Next we establish the connection with the transfinite diameter.
\newline
\newline
{\bf Theorem
\arabic{section}.\arabic{example}.}\addtocounter{example}{1}  {\sl For
any compact set $E\subset \C$,
\begin{equation} \tau(E) = \cp(E).
\label{cap}
 \end{equation}
Moreover, if $E$ has positive capacity,
then
\begin{equation} \nu(\F_n) \stackrel{*}{\rightarrow} \mu_E \
\text{\em  as } n\rightarrow \infty.  \label{positive cap}
\end{equation}
}
\newline
{\bf Proof.}  First we show that
\begin{equation} V_E = \log \frac{1}{\cp(E)} \ge \log \frac{1}{\tau(E)}.
\label{proof eq 1}
 \end{equation} Let $F(z_1, z_2, \ldots, z_n) :=
\sum_{1\le i<j \le n} \log (1/|z_i-z_j|)$.  Then the expected value of
$F$ with respect to the product of equilibrium measures $\dd\mu_E(z_1)
\dd\mu_E(z_2) \cdots \dd\mu_E(z_n)$ cannot be less than its minimum value
defined in (\ref{min prob}); i.e.,
\begin{equation*} \begin{split} \int
\int \cdots \int F(z_1, z_2, \ldots, z_n) \dd\mu_E(z_1)\cdots \dd\mu_E(z_n)
= \frac{n(n-1)}{2} V_E \\ \ge \mathcal{E}_n(E) =  \frac{n(n-1)}{2} \log
\frac {1}{\delta_n(E)}.
  \end{split}
 \end{equation*}
 Dividing by $n(n-1)/2$ and letting $n\rightarrow \infty$ gives
 (\ref{proof eq 1}).

For the reverse direction, let $\hat{\mu}$ be a weak-star limit point of
the measures $\nu_n:=\nu(\F_n)$, say
$\nu_n\stackrel{*}{\rightarrow}\hat{\mu} $ as $n\rightarrow\infty$,
$n\in\mathcal{N}$.  Set $\log_M x:= \min\{\log x, M\}$.  By the Monotone
Convergence Theorem and the weak-star convergence of $\nu_n \times
\nu_n$ to $\hat{\mu}\times\hat{\mu}$ for $n\in\mathcal{N}$, we get
\newline
\newline
$$I(\hat{\mu})  = \int \int \log \frac{1}{|z-t|}
\dd\hat\mu(z)\dd\hat\mu(t) = \lim_{M\rightarrow\infty} \int \int \log_M
\frac{1}{|z-t|} \dd\hat\mu(z) \dd\hat\mu(t) $$
\begin{align*} & =
\lim_{M\rightarrow\infty} \lim_{n\rightarrow\infty} \int \int \log_M
\frac{1}{|z-t|} \dd\nu_n(z) \dd\nu_n(t) \\ & \le \lim_{M\rightarrow\infty}
\lim_{n\rightarrow\infty} \left\{ \frac{2}{n^2} \mathcal{E}_n (E) +
\frac{nM}{n^2} \right\} \\ & = \lim_{M\rightarrow\infty}
\lim_{n\rightarrow\infty} \log \frac{1}{\delta_n(E)} = \log
\frac{1}{\tau(E)} \le V_E.
\end{align*}
Thus from the minimality
property of $V_E$, we have $$V_E \le I(\hat\mu) \le \log
\frac{1}{\tau(E)} \le V_E,$$ which proves (\ref{cap}).  Furthermore, if
$\cp(E)>0$, then by uniqueness of the equilibrium (minimizing) measure,
$\hat\mu = \mu_E$.  Since $\hat\mu$ was an arbitrary limit measure of
$\nu(\F_n)$, (\ref{positive cap}) follows.  $\Box$
\newline

As we have earlier observed, Fekete points necessarily lie on the outer
boundary of $E$.  Thus from (\ref{positive cap}) we immediately deduce
that {\em the equilibrium measure $\mu_E$ is supported on the outer
boundary of $E$}; consequently,
$$\cp(E) = \cp (\partial_\infty E), \ \
\ \ \mu_E = \mu_{\partial_\infty E}.$$
If $\partial_\infty E$ is a
continuum (not a single point), then $\mbox{supp}(\mu_E) =
\partial_{\infty}E$.  In general, $\partial_\infty E \setminus
\mbox{supp}(\mu_E)$ has capacity zero.

From our knowledge of Fekete points for the disk we deduce the
following.
\newline
\newline
{\bf Example \arabic{section}.\arabic{example}.} \addtocounter{example}{1}
If $E$ is the closed disk $|z-a| \le r$, then $\cp(E) = r$ and $\dd\mu_E =
\frac{1}{2\pi r}\dd s$,\footnote{This also follows from the fact that the
disk $E$ is invariant under rotations about $z=a$ and since the
equilibrium measure is unique and supported on the circumference it must
also be rotation invariant and hence of the form described.} where $\dd s$
is arclength on the circumference $|z-a| = r$.  Furthermore, the
equilibrium potential $U^{\mu_E}(z)$ satisfies
\begin{equation}
U^{\mu_E}(z) = \frac{1}{2\pi} \int_0^{2\pi} \log
\frac{1}{|z-a-r\ee^{\ii\theta}|} \dd\theta = \left\{\begin{array}{l l} \log
\frac{1}{r} & \text{ for } |z-a| \le r \\   \log \frac{1}{|z-a|} &
\text{ for } |z-a| \ge r.
  \end{array}\right.  \label{conductor exercise}
 \end{equation}
\begin{figure}[h]
\centering \input{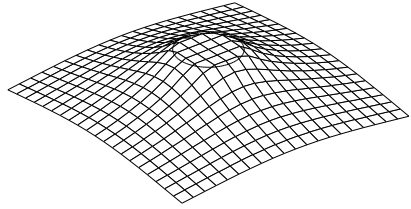}
\caption{Graph of equilibrium potential for the unit disk}
 \end{figure} \label{figeqpot}
\newline
\newline
{\bf \em Exercise.}  Verify
formula (\ref{conductor exercise}).  [Hint: The mean-value property for
harmonic functions is useful here; see Theorem 2.1.]
\newline
\newline
{\bf Example \arabic{section}.\arabic{example}.}
\addtocounter{example}{1} Let $E = [a,b]$ be a segment on the real line.
Then $\cp(E) = (b-a)/4$  and $\dd\mu_E$ is the arcsine measure; i.e.,
$$\dd\mu_E=\frac{1}{\pi}\frac{\dd x}{\sqrt{(x-a)(b-x)}}, \;\;\; x\in[a,b].$$
If  $a = - 1$, $b = 1$, the conductor potential is given by
$$U^{\mu_E}(z)=\log2 \text{ if } x\in[-1,1],$$
$$U^{\mu_E}(z)=\log2-\log\left|z+\sqrt{z^2-1}\right| \text{ if }
z\not\in[-1,1],$$ where the branch $\sqrt{z^2-1}$ is positive for
$z=x>1$.  These facts can be obtained from Example 1.10 by applying the
Joukowski conformal map of $\C\setminus[-1,1]$ onto $|w|>1$.  See
Example 3.6.
\newline

In the examples for the disk and line segment, observe that the
equilibrium potential $U^{\mu_E}$ is constant on $E$, namely it equals
$V_E = \log \frac{1}{\cp(E)}$ there.  This is certainly consistent with
our expectations based on physical grounds, that equilibrium should
occur when the potential (voltage) is constant; for otherwise there
would be a flow of charge to the points of $E$ at lower potential.  From
a mathematically rigorous point of view, this assertion is true
\dword{quasi-everywhere} (\dword{q.e.}) on $E$; that is, except for a set of
capacity zero.\footnote{An arbitrary Borel set $B$ has capacity zero if
$\sup\{\cp(F) : F\subset B \text{ compact}\} = 0$.}  This fact is
included in the following result.
\newline
\newline

{\bf Theorem \arabic{section}.\arabic{example} \addtocounter{example}{1}
(\dword{Frostman's Theorem}).}  {\sl Let $E\subset\C$ be compact with
$\cp(E)>0$.  Then
\begin{enumerate}[(a)]
\item $U^{\mu_E}(z)\le V_E$ for all $z\in\C$;
\item $U^{\mu_E}(z) = V_E$ q.e.\  on $E$.
\end{enumerate}
}

The proof relies on a maximum principle for potentials which is
discussed in Section 2.  For full details, see [R1], [ST], [Ts].

The theorem suggests that we visualize the 3-dimensional graph of an
equilibrium potential as something like an infinite tent with an
``essentially'' flat roof consisting of the projection of the set $E$
and tent sides that flow down and outward to $-\infty$;
see, e.g., Figure \ref{figeqpot}.

There are many important consequences of Frostman's result, of which the
following will be useful in the proof of the main theorem of this
section.
\newline
\newline
{\bf Proposition \arabic{section}.\arabic{example}.} \addtocounter{example}{1}
{\sl Let $E\subset\C$ be compact with $\cp(E)>0$.  If $\sigma$ is any
probability measure with compact support, then}
$$\inf_{z\in E} U^{\sigma}(z)\le V_E = \log \frac{1}{\cp(E)}.$$
\newline
\newline
{\bf Proof.}  Here we use
the \dword{reciprocity law} (a simple consequence of the Fubini-Tonelli
Theorem) which asserts that
$$\int U^\sigma (z)\dd\mu_E(z) = \int U^{\mu_E}(z) \dd\sigma(z).$$
The left-hand side is bounded below by
$\inf_{z\in E} U^\sigma(z)$ and, from Frostman's theorem, $V_E$ is an
upper bound for the right-hand side.  $\Box$
\newline

We now introduce a third quantity associated with a compact set $E$ ---
the \dword{Chebyshev constant}, $\cheb(E)$ --- which arises in a min-max
problem.
\newline
\newline
{\bf Polynomial Extremal Problem:} Determine
the minimal sup norm on $E$ for monic polynomials of degree $n$. That
is, determine\footnote{This problem is equivalent to finding the
polynomial of degree $\le n-1$ of best uniform approximation to the
monomial $z^n$ on $E$.} $$t_n(E):=\min_{p\in\Po_{n-1}}\|z^n+p(z)\|_E,$$
where $\Po_{n-1}$ denotes the collection of all polynomials of degree
$\leq n-1$ and $\|\cdot\|_E$ is the sup norm (uniform norm) on $E$.  We
assume that $E$ contains infinitely many points (which is always the
case if $\cp(E)>0$).  Then for every $n$ there is a unique monic
polynomial $T_n(z)=z^n+\cdots$ such that $\|T_n\|_E=t_n(E)$, which is
called the $n$th \dword{Chebyshev polynomial} for $E$.
\newline
\newline
{\bf \em Exercise.}  Prove that all the zeros of $T_n$ lie in the convex
hull of $E$.  (This fact is due to Fej\'er.)
\newline

In view of the simple chain of inequalities
$$t_{m+n}(E)=\|T_{m+n}\|_E\leq\|T_mT_n\|_E\leq\|T_m\|_E\|T_n\|_E=t_m(E)t_n(E),$$
the sequence $\log t_n (E)$ is subadditive, from which it follows that
$t_n(E)^{1/n}$ converges and its limit is $\inf_{k\ge 1}
\{t_k(E)^{1/k}\}$ (cf.~[Ts], [ST]).  We call this limit the
\dword{Chebyshev constant} for $E$: \begin{equation} \cheb(E) :=
\lim_{n\rightarrow\infty} t_n(E)^{1/n} = \inf_{k\ge 1}\{t_k(E)^{1/k}\}.
\label{cheb}
 \end{equation}

From (\ref{cheb}) and the definition of $t_n(E)$ we deduce the
following.
\newline
\newline
 {\bf Lemma
\arabic{section}.\arabic{example}.} \addtocounter{example}{1}  {\em For
any monic polynomial $p_n(z)$ of degree $n$ there holds }
\begin{equation} \lVert p_n\rVert_E \ge [\cheb(E)]^n.  \label{monic
cheb}
 \end{equation}
\newline
 {\bf Example
\arabic{section}.\arabic{example}.} \addtocounter{example}{1} Let $E$ be
the closed disk of radius $R$, centered at 0.  For any $p\in\Po_{n-1}$,
the ratio $(z^n+p(z))/z^n$ represents an analytic function in $|z|\geq1$
that takes the value 1 at $\infty$. By the maximum principle for
analytic functions,
$$\|z^n+p(z)\|_E=\max_{|z|=R}|z^n+p(z)|=R^n\max_{|z|=R}\left|\frac{z^n+p(z)}{z^n}\right|\geq
R^n,$$ and strict inequality takes place if $p(z)$ is not identically
zero.  It follows that $T_n(z)=z^n$. Therefore $t_n(E)=R^n$ and
$\cheb(E)=R$.
\newline
\newline
 {\bf Example
\arabic{section}.\arabic{example}.} \addtocounter{example}{1}  Let $E=
[-1, 1]$. Then $T_n$ is the classical monic Chebyshev
polynomial\footnote{There is ambiguity in this notation since $T_n(x)$
is traditionally used to denote $\cos(n\arccos x)$.}
$$T_n(x)=2^{1-n}\cos(n\arccos x), \;\;\; x\in[-1,1], \;\;\; n\geq1.$$
Thus, $t_n(E)=\lVert T_n \rVert_{[-1,1]} = 2^{1-n}$ from which it
follows that $\cheb(E)=1/2$, which is the same as the capacity of $[-1,
1]$.
\newline
\newline
{\bf \em Exercise.}  Verify
Example 1.16 by using the fact that $2^{1-n}\cos(n\arccos x)$
equioscillates $n+1$ times on $[-1,1]$.
\newline

Closely related to Chebyshev polynomials are {\em Fekete polynomials}.
An $n$th degree Fekete polynomial $F_n(z)$ is a monic polynomial having
all its zeros at the $n$ points of a Fekete set $\F_n$.
\newline
\newline

 {\bf Example \arabic{section}.\arabic{example}.}
\addtocounter{example}{1} If $E$ is the closed unit disk centered at 0,
then one can take $F_n(z) = z^n-1$, so that $\|F_n\|_E=2$. Comparing
this with Example 1.15
we see that the $F_n$'s are
asymptotically optimal for the Chebyshev problem:
$$\lim_{n\to\infty}\|F_n\|_E^{1/n}=\lim_{n\to\infty}\|T_n\|_E^{1/n}=1=\cheb(E).$$
Moreover, uniformly on compact subsets of $|z| > 1$, we have
$$\lim_{n\to\infty}|F_n(z)|^{1/n}=\lim_{n\to\infty}|T_n(z)|^{1/n}=|z|=\exp(-U^{\mu_E}(z)),$$
(the last equality follows from formula (\ref{conductor exercise})).
\newline

The above examples illustrate the following fundamental theorem, various
parts of which are due to Fekete, Frostman, and Szeg\H{o}.
\newline

{\bf Theorem \arabic{section}.\arabic{example}}
\addtocounter{example}{1} (\dword{Fundamental Theorem of Classical Potential
Theory}).  {\sl For any compact set $E\subset\C$,
\newline
\textnormal{(a)}  $\cp(E)=\tau(E)=\cheb(E)$;
\newline
 \textnormal{(b)}
Fekete polynomials are asymptotically optimal for the Chebyshev problem:
$$\lim_{n\to\infty}\|F_n\|_E^{1/n}=\cheb(E)=\cp(E).$$ If $\cp(E ) > 0$
(so that $\mu_E$ is unique), then we also have\textnormal{:}
\newline
\textnormal{(c)} Fekete points (the zeros of  $F_n$) have asymptotic
distribution $\mu_E$, i.e., $\nu(\F_n) \stackrel{*}{\rightarrow} \mu_E$
as $n\rightarrow\infty$;
\newline
 \textnormal{(d)} uniformly on compact
subsets of the unbounded component of  $\C \setminus E$,
$$\lim_{n\to\infty}|F_n(z)|^{1/n}=\exp(-U^{\mu_E}(z)).$$}
\newline
\newline
{\bf Proof.}  Part (c) was established in
(\ref{positive cap}) of Theorem 1.9.  Assertion (d) follows from (c) on
observing that $$\frac{1}{n}\log\frac{1}{|F_n(z)|} = U^{\nu(\F_n)}(z),$$
and that all the Fekete points lie on $E$.  Regarding (a) and (b), we
already know that $\cp(E) = \tau(E)$.  So to establish (a), we prove that
$\tau(E) = \cheb(E)$.

Let $\F_n = \{z_{k}^{(n)}\}_{k=1}^{n}$ denote an $n$-point Fekete set
for $E$.  Then $$\delta_{n+1}^{n(n+1)/2} = \max_{\{z_i\}\subset E}
\prod_{1\le i<j\le n+1} |z_i-z_j| \ge \left[\prod_{k=1}^n
|z-z_k^{(n)}|\right] \delta_{n}^{n(n-1)/2}$$ for all $z\in E$.  Thus
$$\delta_{n+1}^{n(n+1)/2}/\delta_{n}^{n(n-1)/2} \ge |F_n(z)|,\ \ \ \ \
z\in E,$$ and so on taking $n$th roots, we get \begin{equation}
\delta_{n+1}^{(n+1)/2}/\delta_{n}^{(n-1)/2} \ge \lVert F_n
\rVert_E^{1/n}.  \label{n-th roots delta}
 \end{equation} The left-hand
side of this inequality can be written as
$$\left(\frac{\delta_{n+1}}{\delta_n}\right)^{n/2}
\delta_{n+1}^{1/2}\delta_{n}^{1/2},$$ which is bounded above by
$\delta_{n+1}^{1/2}\delta_{n}^{1/2}$ since the sequence $\delta_n$ is
decreasing.  From (\ref{monic cheb}), we have that the right-hand side of
(\ref{n-th roots delta}) is bounded below by $\cheb(E)$.  Hence
$$\delta_{n+1}^{1/2}\delta_{n}^{1/2} \ge  \lVert F_n \rVert_E^{1/n} \ge
\cheb(E),$$ and so on letting $n\rightarrow\infty$, we get $$\tau(E)\ge
\limsup_{n\rightarrow\infty}\, \rVert F_n \lVert_E^{1/n}\, \ge
\liminf_{n\rightarrow\infty}\, \rVert F_n \lVert_E^{1/n}\,
\ge\cheb(E).$$

It remains only to prove that $\cheb(E) \ge \tau(E)$.  This is obvious
if $\tau(E)=\cp(E)=0$.  So assume $\cp(E)>0$.  Let $\nu(T_n)$ denote the
normalized counting measure in the zeros of $T_n(z)$.  Then by
Proposition 1.13,
$$\inf_{z\in E} U^{\nu(T_n)}(z) =
\inf_{z\in E} \frac{1}{n}\log\frac{1}{|T_n(z)|} =
\frac{1}{n}\log\frac{1}{t_n(E)}\le V_E.$$ Hence $$t_n(E)^{1/n} \ge
\ee^{-V_E} = \cp(E) = \tau(E),$$ and on letting $n\rightarrow\infty$, we
get $\cheb(E) \ge \tau(E)$.  $\Box$
\newline
\newline
 {\bf \em
Exercise.}  Let $0<a<b$.  Prove that $$\cp\left([a,b]\cup[-b,-a]\right)
= \frac{\sqrt{b^2-a^2}}{2}.$$ [Hint: What are the Chebyshev polynomials
of even degree for this union?]
\newline
\newline
 {\bf \em Exercise.}
Let $P(z)$ be a monic polynomial of degree $n$, and consider the {\em
lemniscate set} $L:= \{z : |P(z)|\le R^n\}$.  Prove that $\cp(L) = R$.
[Hint: Begin by determining the Chebyshev polynomials $T_{kn},\ k=1, 2,
\ldots, $ for $L$.]
\newline
\newline
 {\bf \em Exercise.}  Let $E$ be a
compact set and $\epsilon>0$.  Show that there exists a lemniscate set
$L$ such that $E\subset L$ and $\cp(L)<\cp(E)+\epsilon$.

\sect{Harmonic, Superharmonic and Subharmonic Functions}
\setcounter{example}{1} \setcounter{equation}{0}

Recall that a real-valued function $u(z)$ defined in an open set
$D\subset\C$ is \dword{harmonic} in $D$ if $u$ and its 1st and 2nd partial
derivatives are continuous in $D$ and $u$ satisfies {\em Laplace's
equation} \begin{equation} u_{xx}(z) + u_{yy}(z) = 0,\ \ \ \ \ z\in D.
\label{laplace}
 \end{equation} (Actually it is enough to merely assume
that the 2nd partial derivatives exist and satisfy (\ref{laplace}).)
{\em Locally}, harmonic functions are the real (or imaginary) parts of
an analytic function.
\newline
\newline
 {\bf \em Exercise.}  Prove that
if $u$ is harmonic in $D$, then $g(z) := u_x(z) - iu_y(z)$ is analytic
in $D$.
\newline

Note that the important function $\log|z|$ is harmonic in $\C\setminus
\{0\}$, is locally the real part of a branch of $\log z$, but is not
globally (in $\C\setminus \{0\}$) the real part of an analytic function.
\newline
\newline
 {\bf \em Exercise.}  Prove that the logarithmic
potential $U^{\mu}(z)$ is harmonic for $z$ not in the support of $\mu$
(assuming this support is compact and $\mu$ is a finite measure).
\newline

Harmonic functions can also be characterized by the following
\dword{Mean-Value Property} (\dword{MVP}).
\newline
\newline
 {\bf Theorem
\arabic{section}.\arabic{example}.} \addtocounter{example}{1}  {\sl A
real-valued function $u(z)$ is harmonic in an open set $D$ if and only
if $u$ is continuous in $D$ and locally satisfies the \textbf{mean-value
property}; i.e., if the disk $|z-a|\le r$ is contained in $D$, then}
\begin{equation} u(a) = \frac{1}{2\pi} \int_0^{2\pi} u(a+r\ee^{\ii\theta})
\dd\theta.  \label{mvp}
 \end{equation} (In fact, it is enough that
equality holds for $r=r(a)$ sufficiently small.)
\newline
\newline
 {\bf
\em Exercise.}  Prove that if $u$ is harmonic on an open set $D$
containing the disk $|z-z_0|\le r$, then $$u(z_0) = \frac{1}{\pi r^2}
\int \int_{|z-z_0|\le r} u(z) \dd x \dd y.$$
\newline
\newline
 {\bf \em
Exercise.}  Prove that if $u_n(z)$, $n=1,2,\ldots$, is a sequence of
functions harmonic in $D$ that converges locally uniformly to a function
$u$ in $D$, then $u$ is harmonic in $D$.
\newline

An important consequence of (\ref{mvp}) is the \dword{Max-Min Principle
for harmonic functions}.
\newline
\newline
 {\bf Theorem \arabic{section}.\arabic{example}.} \addtocounter{example}{1}
 {\sl If $u$ is harmonic in a \textbf{domain} $D$ (i.e.,  an open connected set)
and $u$ attains its maximum (or minimum) in $D$, then $u$ is identically
constant in $D$.

Furthermore, if $u$ is harmonic in the interior and continuous on the
boundary of a compact set, then $u$ attains its max and min on the
boundary.}
\newline

The proof of this principle follows from the observation that if $u$
attains its max at a point $z_0\in D$ and $\overline{D}_r(z_0)$ is any
closed disk centered at $z_0$, then $u$ must equal $u(z_0)$ for all
$z\in\overline{D}_r(z_0)$ since the contrary assumption would lead to a
violation of the MVP (simply integrate around a circle centered at $z_0$
and containing a point where $u(z)<u(z_0)$).

Theorem 2.2 also tells us that harmonic functions
are determined by their values on the boundary of a compact set.
Indeed, in the case of a disk, we have \dword{Poisson's integral formula}:
If $u$ is harmonic in $|z|<R$ and continuous on $|z|\le R$, then
\begin{equation} u(z) = \frac{1}{2\pi} \int_0^{2\pi} P(t,z)u(t)\dd\theta,
\ \ \ \ \ t=R\ee^{\ii\theta}, \ \ \ \ |z|<R, \label{poisson 1}
\end{equation} where \begin{equation} P(t,z) :=
\frac{|t|^2-|z|^2}{|t-z|^2} = \Re\left(\frac{t+z}{t-z}\right).
\label{poisson 2}
 \end{equation}

This formula can be deduced, for example, from the Cauchy integral
formula for analytic functions; it includes as a special case the MVP
(\ref{mvp}).
\newline
\newline
 {\bf \em Exercise.}  Prove that if
$U(t)$ is integrable (in the Lebesgue sense) on $|t|=R$, then
\begin{equation} u(z) := \frac{1}{2\pi} \int_0^{2\pi} P(t,z)U(t)\dd t
\label{integrable}
 \end{equation} is harmonic in $|z|<R$.  (If $U(t)$ is
continuous on the circle $|t|=R$, then Schwarz's theorem asserts that
$u$ as given in (\ref{integrable}) solves the \emph{Dirichlet problem}
for the disk; i.e.,  $\lim_{z\to t}u(z)=U(t)$ for all $t$ on the boundary
$|t|=R$.)
\newline

If we replace equality in (\ref{mvp}) by $\ge$, we obtain the class of
superharmonic functions.
\newline
\newline
 {\bf Definition
\arabic{section}.\arabic{example}.} \addtocounter{example}{1} An
extended real-valued function $f$ on an open set $D\subset\C$ is
called\linebreak \dword{superharmonic} in $D$ if $f$ is not the constant
function $+\infty$ and satisfies
\begin{enumerate}[(i)]
\item $f$ is
lower semi-continuous on $D$;
\item the value of $f$ at any point
$z_0\in D$ is not less than its average over any circle in $D$ centered
at $z_0$; that is \begin{equation} f(z_0) \ge \frac{1}{2\pi}
\int_0^{2\pi} f(z_0 + r\ee^{\ii\theta})\dd\theta \label{superharmonic def}
\end{equation} provided the closed disk $\overline{D_r(z_0)} := \{z\in\C
: |z-z_0|\le r\}$ is contained in $D$.
  \end{enumerate}

A \dword{subharmonic} function is the negative of a superharmonic
function.  A real-valued function $u(z)$ that is both superharmonic and
subharmonic in $D$ is harmonic in $D$.
\newline
\newline
 {\bf \em
Exercise.}  Let $f:D\rightarrow\R$, $D$ a domain, $f\in C^2(D)$.  Prove
that $f$ is superharmonic in $D$ if and only if $\Delta f := f_{xx} +
f_{yy} \le 0$ at all points of $D$.  [Hint: Begin by showing that a
negative Laplacian implies that $f$ is superharmonic.]
\newline

As the above exercise illustrates, superharmonic functions in $\R^2$ are
analogues of concave functions in $\R$; and subharmonic functions are
the analogues of convex functions.
\newline
\newline
 {\bf \em
Exercise.}  Prove that if $F(z)$ is analytic in a domain $D$ and $p>0$,
then $|F(z)|^p$ is subharmonic in $D$.  Furthermore, show $\log|F(z)|$
is subharmonic in $D$ unless $F$ is identically zero.
\newline

Essential for us is the fact that {\sl logarithmic potentials are
superharmonic in} $\C$.  Indeed, as we have seen, $$U^\mu (z) = \int
\log \frac{1}{|z-t|}\dd\mu(t)$$ is l.s.c.\  on $\C$ and since $\log
(1/|z-t|)$ is superharmonic in $\C$ for fixed $t$ (Why?), it follows
from the Fubini-Tonelli theorem that for $a\in\C$ \begin{align*}
\frac{1}{2\pi}\int_0^{2\pi} U^\mu(a+r\ee^{\ii\theta})\dd\theta & = \int
\frac{1}{2\pi}\int_{-\pi}^{\pi} \log \frac{1}{|a+r\ee^{\ii\theta} -t|}
\dd\theta \dd\mu(t) \\ & \le \int \log \frac{1}{|a-t|}\dd\mu(t) = U^\mu(a)
\end{align*} (see also (\ref{conductor exercise})).
\newline
\newline
{\bf \em Exercise.}  Prove that if $U^\mu$ is harmonic in a neighborhood
of a point $z_0$, then $z_0\not\in\mbox{supp}(\mu)$.
\newline

While it may appear that potentials are rather special types of
superharmonic functions, their properties are key to the analysis of
general superharmonic functions.  This is thanks to the following
celebrated result.
\newline
\newline
 {\bf Theorem \arabic{section}.\arabic{example}} \addtocounter{example}{1}
 (\dword{Riesz Decomposition}).  {\sl If $f$ is superharmonic in a domain $D$,
 then there exists a positive measure $\lambda$ supported on $D$ such that for
every subdomain $D^*\subset D$  for which $\overline {D^*}\subset D$, we
have \begin{equation} f(z) = h(z) + \int_D \log
\frac{1}{|z-t|}\dd\lambda(t), \ \ \ z\in D^*, \label{riesz theorem}
\end{equation} where $h$ is harmonic in $D^*$.}
\newline

For the case when $f$ is smooth, superharmonicity implies that $\Delta f
= f_{xx} + f_{yy} \le 0$ in $D$ and it turns out that the positive
measure \begin{equation} \boxed{\dd\lambda(t) := -\frac{1}{2\pi}\Delta
f(t)\dd m_2(t), \ \ \ \ \ t\in D,} \label{riesz smooth}
 \end{equation}
where $m_2$ denotes 2-dimensional Lebesgue measure, yields the
appropriate potential $U^\lambda$ for which (\ref{riesz theorem}) holds.

Let's verify this for the simple but important case when $f(z) =
-|z|^2$, for which we have $\Delta f \equiv -4$.  Then we need to show
that $\Delta (f-U^\lambda) = 0$ in $D^*$ where $\dd\lambda = (2/\pi)\dd m_2$.
For an arbitrary disk $D_r(z_0):=\{z: |z-z_0| < r \}$ contained in
$D^*$, write $U^\lambda = U^{\lambda_1}+U^{\lambda_2}$ where $\lambda_1
= \lambda|_{D_r(z_0)}$ and $\lambda_2 = \lambda - \lambda_1$.  Then
$U^{\lambda_2}$ is harmonic in $D_r(z_0)$ and so we need only show that
$\Delta(f-U^{\lambda_1}) = 0$ in $D_r(z_0)$.  A simple calculation using
(\ref{conductor exercise}) gives that 
\begin{align*}
U^{\lambda_1}(z) & = \frac{2}{\pi}\int_{D_r(z_0)} \log
\frac{1}{|z-t|}\dd m_2(t) \\ & = 2r^2 \log \frac{1}{r} + r^2 - |z-z_0|^2,
\end{align*}
from which we get $\Delta U^{\lambda_1}(z) = -4 = \Delta f$
for $z\in D_r(z_0)$, as desired.   For more general but smooth $f$, one
can use Green's formula\footnote{Recall that Green's formula for smooth
functions $u, v$ on a bounded open set $D$ with smooth boundary
$\partial D$ asserts that $$\int_D (v\Delta u - u\Delta v)\dd m_2 =
-\int_{\partial D} (v\frac{\partial u}{\partial n} - u\frac{\partial
v}{\partial n})\dd s,$$ where $\partial/\partial n$ denotes
differentiation in the direction of the inner normal to $D$.}  to verify
that (\ref{riesz smooth}) yields the decomposition in (\ref{riesz
theorem}).

For general superharmonic functions $f$, we interpret the right-hand side
of (\ref{riesz smooth}) in the {\em distributional sense} (cf.~[R1,
Sec.~3.7]); more precisely, we identify $-\Delta f\dd m_2$ as the unique
positive measure that satisfies \begin{equation} \int_D \phi(-\Delta
f\dd m_2) := -\int_D f\Delta \phi \dd m_2 \label{neg Delta f}
\end{equation} for all $C^\infty$ functions $\phi$ whose support is a
compact subset of $D$.  This condition is precisely what would be
expected from Green's formula.  The existence of such a measure $-\Delta
f\dd m_2$ satisfying (\ref{neg Delta f}) is guaranteed by the Riesz
representation theorem for linear functionals.  With this
interpretation, it follows that for any finite Borel measure $\mu$ with
compact support there holds
\begin{equation} \boxed{\mu =
-\frac{1}{2\pi}\Delta U^\mu.} \label{Borel mu}
\end{equation}

Just as the MVP (\ref{mvp}) for harmonic functions
implied the Max-Min Principle (Theorem 2.2), the mean-value inequality
property (\ref{superharmonic def}) yields the following.
\newline
\newline
 {\bf Theorem \arabic{section}.\arabic{example}}
\addtocounter{example}{1} (\dword{Minimum Principle for Superharmonic
Functions}).  {\sl Let $D$ be a bounded domain and $g$ a superharmonic
function on $D$ such that
\begin{equation} \liminf_{z\rightarrow\zeta}
g(z) \ge m \ \ \ \ \forall \ \zeta \in \partial D.
\label{min subharmonic}
\end{equation}
Then $g(z)>m$ for all $z$ in $D$, unless $g$ is constant.}
\newline
\newline
 {\bf \em Exercise.}  Use the Min
Principle to prove that a l.s.c.\  function $f$ (not identically
$+\infty$) is superharmonic in a domain $D$ if and only if it has the
following property:  If $D_0\subset D$ is a bounded domain whose closure
is contained in $D$ and $u$ is harmonic in $D_0$, continuous on
$\overline{D_0}$ and satisfies $u(\zeta)\le f(\zeta)\ \ \forall\ \zeta
\in \partial D_0$, then $u(z) \le f(z) \ \ \forall \ z\in D_0$.
\newline

A more general form of the Min Principle allows us to ignore a set of
points on $\partial D$ of capacity zero provided $g$ is lower bounded on
$D$; see [ST].
\newline
\newline
 {\bf Theorem \arabic{section}.\arabic{example}} \addtocounter{example}{1}
(\dword{Generalized Min Principle}).  {\sl If $D\subset\overline{\C}$ is a
domain, $\cp(\partial D)>0$, $g$ is superharmonic and bounded from below
in $D$ and (\ref{min subharmonic}) holds for q.e.\  $\zeta$ on
$\partial D$, then $g(z)>m\ \ \forall \ z\in D$ unless $g$ is constant.}
\newline
\newline
 {\bf \em Exercise.} Give an example to show that the
lower boundedness assumption cannot be removed in the above result.
\newline

For potentials what is crucial is their behavior on the support of its
defining measure.
\newline
\newline
 {\bf Theorem \arabic{section}.\arabic{example}}\addtocounter{example}{1}
(\dword{Maximum Principle for Potentials}).  {\sl Let $\mu$ be a finite
positive measure with compact support.  If $U^\mu (z) \le M$ for all
$z\in\textnormal{supp}(\mu)$, then $U^\mu (z) \le M$ for all $z\in\C$.}
\newline

The Max Principle is a special case of the following important result.
\newline
\newline
 {\bf Theorem \arabic{section}.\arabic{example}}
\addtocounter{example}{1} (\dword{Principle of Domination}).
{\sl Let $\mu$, $\nu$ be positive finite measures with compact supports,
$\nu(\C)\le\mu(\C)$, and $\mu$ has finite logarithmic energy
$(I(\mu)<\infty)$.  If for some constant $c$, the inequality $$U^\mu(z)
\le U^\nu (z) + c$$ holds $\mu$-a.e., then it holds for all $z\in\C$.}
\newline

The idea of the proof of the above theorem is to consider the function
$$U(z):= \min (U^\nu (z) +c, U^\mu (z)),$$ which is superharmonic since
the minimum of two superharmonic functions is again superharmonic.  Let
$\lambda:= -\frac{1}{2\pi} \Delta U$ (i.e., $\lambda$ is the measure
guaranteed by the Riesz Decomposition Theorem (RDT) with $D=\C$) and
argue that $\lambda$ must equal $\mu$.  See [ST] for details.

As an application of the Principle of Domination, we present
\newline
\newline
 {\bf Example \arabic{section}.\arabic{example}.}
\addtocounter{example}{1}  Let $p_n$, $n=1, 2, \ldots$, be a sequence of
monic polynomials of respective degrees $n$ that satisfy
\begin{equation} \limsup_{n\rightarrow\infty}\,\lVert
p_n\rVert_{[-1,1]}^{1/n}\, \le \frac{1}{2} = \cp([-1,1]), \label{monic
poly seq}
 \end{equation} and let $\nu(p_n)$ denote the normalized
counting measure in the zeros of $p_n$.  If all the zeros of the $p_n$'s
lie on $[-1,1]$, then
\begin{equation}
\nu(p_n)\stackrel{*}{\rightarrow}\mu_{[-1,1]} =
\frac{\dd x}{\pi\sqrt{1-x^2}}\ \ \text{ as }\  n\rightarrow\infty.
\label{nu limit}
\end{equation}
Indeed, by definition of the sup norm,
$$\frac{1}{n}\log\frac{1}{|p_n(z)|}\ge\frac{1}{n}\log\frac{1}{\lVert
p_n\rVert_{[-1,1]}},\ \ \ z\in[-1,1].$$ We write this last inequality in
the equivalent form \begin{equation} U^{\nu(p_n)}(z) + \log2 \ge
\frac{1}{n}\log \frac{1}{\lVert p_n\rVert_{[-1,1]}} +
U^{\mu_{[-1,1]}}(z),\ \ \ z\in[-1,1], \label{equiv nu}
 \end{equation}
where  we used the fact that $U^{\mu_{[-1,1]}}(x) =
\log2$  for all $ x\in[-1,1]$ (see Examples 1.11 and 3.6).  By the
Principle of Domination, (\ref{equiv nu}) holds for all $z\in\C$.  Now
let $\nu$ be a weak-star limit measure of the sequence $\nu(p_n)$.  Then
from (\ref{monic poly seq}) and (\ref{equiv nu}), we get $$U^\nu(z) +
\log 2 \ge \log2 + U^{\mu_{[-1,1]}}(z) \ \ \ \forall\, z\not\in[-1,1],$$
i.e., $U^\nu(z) - U^{\mu_{[-1,1]}}(z)\ge0$ for all
$z\in\Omega:=\overline{\C}\setminus[-1,1]$.  But since $U^\nu -
U^{\mu_{[-1,1]}}$ vanishes at $\infty$ and is harmonic in $\Omega$, the
Max-Min Principle asserts that $U^\nu(z) = U^{\mu_{[-1,1]}}(z)$ for
$z\in\Omega$.  By l.s.c.,  we then deduce that $U^\nu(x) \le
U^{\mu_{[-1,1]}}(x) = \log2$  for all $ x\in[-1,1]$ and so
$I(\nu)\le\log2 = I(\mu_{[-1,1]})$.  By uniqueness of the minimizing
measure, we get that $\nu = \mu_{[-1,1]}$, which proves (\ref{nu limit}).

Actually, (\ref{nu limit}) {\sl holds without any prior assumptions on
the location of zeros of the $p_n$'s.}  This follows from the fact that
(\ref{monic poly seq}) implies that the proportion of zeros that lie
outside any open set containing $[-1,1]$ is asymptotically negligible
(see Section 3).
\newline

Another consequence of the Riesz Decomposition Theorem is the following.
\newline
\newline
 {\bf Theorem \arabic{section}.\arabic{example}}
\addtocounter{example}{1} (\dword{Unicity Theorem}).  {\sl Let $\mu$, $\nu$ be
positive finite measures having compact support.  If, in a region
$D\subset\C$, there holds $$U^\mu (z) = U^\nu (z)+ h(z) \ \ \ \
m_2\mbox{-a.e.,}$$ where $h$ is harmonic in $D$, then $\mu|_D =
\nu|_D$.}
\newline

In particular, if two potentials $U^\mu$ and $U^\nu$ agree except for a
set of 2-dimensional Lebesgue measure zero, then $\mu=\nu$.

\sect{Equilibrium Potentials, Green Functions and Regularity}
\setcounter{example}{1} \setcounter{equation}{0}

Throughout this section, $E\subset\C$ denotes a compact set with
$\cp(E)>0$.  Here, we discuss properties and characterizations of the
equilibrium (conductor) potential $U^{\mu_E}$.

According to Frostman's Theorem 1.12, $U^{\mu_E}$ is
``essentially'' constant on $E$ (more precisely, constant {\em q.e.} on
$E$).  So the following result should come as no surprise.
\newline
\newline
 {\bf Theorem \arabic{section}.\arabic{example}.}
\addtocounter{example}{1}  {\sl If $\nu\in\M(E)$ has finite logarithmic
energy {\rm(i.e.,  $I(\nu)<\infty$)} and $U^\nu(z) = c$ q.e.\  on $E$,
then $c=V_E$ and $\nu= \mu_E$.}
\newline

In the proof of this result, a useful fact is the following.
\newline
\newline
 {\bf \em Exercise.}  If a measure $\nu$ has finite logarithmic
energy, then any set of capacity zero has $\nu$ measure zero.
\newline

With this fact, Theorem 3.1 follows by simply
integrating the equality $U^\nu=c$ with respect to $\dd\mu_E$ and
interchanging order of integration.
\newline
\newline
 {\bf \em
Exercise.}  Give an example to show that if $\nu\in\M(E)$ and $U^\nu$ is
constant on $E$ {\em except for one point}, then $\nu$ need not equal
$\mu_E$.
\newline

What can be said about the continuity properties of $U^{\mu_E}$?
Certainly $U^{\mu_E}$ is continuous in $\C\setminus\mbox{supp}(\mu_E)$
since it is harmonic there.  So suppose $z_0\in\mbox{supp}(\mu_E)$.
Then by l.s.c.\ and Frostman's Theorem, we have $$U^{\mu_E}(z_0) \le
\liminf_{z\rightarrow z_0} U^{\mu_E}(z) \le \limsup_{z\rightarrow z_0}
U^{\mu_E}(z) \le V_E.$$ Hence if $U^{\mu_E}(z_0) = V_E$, then
$U^{\mu_E}$ is continuous at $z_0$.  The converse is also true.  (Prove
it!)  To summarize, we have:
\newline
\newline
 {\bf Theorem \arabic{section}.\arabic{example}.} \addtocounter{example}{1}
{\sl $U^{\mu_E}$ is continuous at $z_0\in\textnormal{supp}(\mu_E)$ if and
only if $U^{\mu_E} (z_0)= V_E$.  Consequently, $U^{\mu_E}$ is continuous
q.e.\  in the plane.}
\newline
\newline
 {\bf Definition
\arabic{section}.\arabic{example}.} \addtocounter{example}{1}  A point
$z_0\in\partial_\infty E$ is said to be a \dword{regular point} of the
unbounded component $\Omega$ of $\overline{\C}\setminus E$ if
$U^{\mu_E}(z_0) = V_E$.  Otherwise, $z_0$ is called an \dword{irregular
point}.  (From Theorem 3.2, we see that the set of all irregular points
has capacity zero.)  If every point of $\partial\Omega = \partial_\infty
E$ is regular, we say that $\Omega$ is regular (with respect to the
Dirichlet problem).
\newline
\newline
 {\bf \em
Exercise.}  Prove that every interior point of $E$ satisfies
$U^{\mu_E}(z) = V_E$.
\newline

The equilibrium potential is related to the Green function associated
with the unbounded component of the complement of $E$; more precisely,
we have
\newline
\newline
 {\bf Definition \arabic{section}.\arabic{example}.} \addtocounter{example}{1}
The \dword{Green function with pole at $\infty$} for the unbounded component
$\Omega$ of $\C\setminus E$ is defined by
\begin{equation} g_\Omega (z, \infty) := V_E - U^{\mu_E}(z).  \label{green infinity}
\end{equation}
(Some authors write $g_E (z, \infty)$ instead of $g_\Omega (z, \infty).
$)
\newline

Three properties uniquely {\em characterize} this function for
$z\in\Omega$; namely
\begin{enumerate}[(a)]
\item $g_\Omega (z, \infty)$ is harmonic in $\Omega \setminus \{\infty\}$ and
bounded from above and below outside each neighborhood of $\infty$;
\item $g_\Omega (z, \infty) - \log |z| = O(1)$ as $z\rightarrow\infty$;
\item $g_\Omega (z, \infty)\rightarrow0$ as $z\rightarrow\zeta, \ z\in \Omega$,
for q.e.\  $\zeta\in\partial\Omega$.
\end{enumerate}

Regarding property (b), it follows from (\ref{green infinity}) and the
fact that $\mu_E$ is a unit measure that \begin{equation}
g_\Omega(z, \infty) - \log|z| \rightarrow V_E = \log
\frac{1}{\cp(E)}\mbox{ as } z\rightarrow\infty.  \label{green limit}
\end{equation} It is also clear from (\ref{green infinity}) that
$g_\Omega (z, \infty) \ge 0$, $g_\Omega (z,\infty) > 0$ for $z\in\Omega$,
and in view of Theorem 3.2, if $\zeta\in\partial\Omega = \partial_\infty
E$, then $g_\Omega(z, \infty) \rightarrow 0$ as $z\rightarrow\zeta$,
$z\in\Omega$, if and only if $\zeta$ is a regular point of $\Omega$.
\newline
\newline
 {\bf \em Exercise.}  Prove that if
$\tilde{g}(z)$ is a function that satisfies properties (a), (b), and (c)
for $z\in\Omega$, then $\tilde{g}(z) = g_\Omega(z, \infty)$ for
$z\in\Omega$.
\newline

In the case when $\Omega$ is {\em simply connected}, we can relate the
Green function to a Riemann mapping of $\Omega$.
\newline
\newline
{\bf Theorem \arabic{section}.\arabic{example}.} \addtocounter{example}{1}
{\sl If the unbounded component $\Omega$ of $\C\setminus E$ is simply
connected, then $g_\Omega (z, \infty) = \log |\Phi(z)|, \ z\in\Omega$,
where $w=\Phi(z)$ is the unique Riemann mapping function from $\Omega$
to the exterior $\{w : |w| > 1\}$ of the unit disk  such that
$\Phi(\infty) = \infty, \ \Phi'(\infty)>0$.}
\newline

Such a function $\Phi$ has a Laurent expansion about $\infty$ of the
form \begin{equation} \Phi(z) = \frac{z}{c} + a_0 + \frac{a_1}{z} +
\frac{a_2}{z^2} + \cdots, \text{ with }c>0, \label{Phi laurent}
\end{equation} and using this representation, properties (a), (b), and
(c) are easy to establish.
\newline
\newline
 {\bf \em Exercise.}  Prove Theorem 3.5.
\newline

From (\ref{Phi laurent}), we immediately see that $$\log |\Phi(z)| -
\log|z| \rightarrow \log\frac{1}{c} \text{ as } z\rightarrow \infty,$$
and comparison with (\ref{green limit}) shows that $$\boxed{c=\cp(E).}$$
From this fact, we can determine the capacity of any compact set $E$
providing we know the exterior conformal mapping function $\Phi$.  We
illustrate this for the line segment.
\newline
\newline
 {\bf Example
\arabic{section}.\arabic{example}.} \addtocounter{example}{1}  Let $E  =
[-1,1]$.  The well-known Joukowski transformation $z = \psi(w) =
\frac{1}{2}(w+w^{-1})$ maps the exterior of the unit circle onto $\Omega
:= \overline{\C} \setminus [-1,1]$, with $\psi(\infty) = \infty, \
\psi'(\infty)>0$.  Solving for $w$, we obtain the desired Riemann
mapping: $$w=\Phi(z) = z+\sqrt{z^2-1},$$ where $\sqrt{z^2-1}$ behaves
like $z$ near infinity.  Thus $$g_\Omega(z,\infty) =
\log\left|z+\sqrt{z^2-1}\right|,$$ and since $\Phi(z) = 2z + \cdots$
near infinity, we get that $\cp([-1,1]) = 1/2$.  Furthermore, from
(\ref{green infinity}), we have $$U^{\mu_E}(z) =
\log2-\log\left|z+\sqrt{z^2-1}\right|$$ as claimed in Example 1.11.
\newline
\newline
 {\bf \em Exercise.}  Show that if
the unbounded component $\Omega$ of $\overline{\C}\setminus E$ is simply
connected, then every point of $\partial_\infty E$ is regular, i.e.,
$\Omega$ is a regular domain.
\newline
\newline
 {\bf \em Exercise.} By
constructing a suitable mapping function, show that the capacity of an
ellipse with semi-axis lengths $a$ and $b$ is $(a+b)/2$.
\newline
\newline
 {\bf \em Exercise.} Show that if the compact set $E$ has
positive capacity and $p(z)$ is a monic polynomial of degree $n$, then
the set $p^{-1}(E)$ has capacity $[\cp(E)]^{1/n}$.
\newline
\newline

In the case when $\partial_\infty E$ is a smooth closed Jordan curve,
there is a simple representation for $\mu_E$ in terms of the Green
function $g_\Omega = g_\Omega(z,\infty)$. Using Green's formula, one
first shows that, for $z \in \Omega$, the equilibrium potential
($=V_E-g_{\Omega}$) identically equals the potential of the unit measure
$\frac{1}{2\pi}\frac{\partial g_\Omega}{\partial n}\dd s$, where the
derivative is taken in the direction of the {\em outer} normal on
$\partial_\infty E $ and $s$ denotes arclength on $\partial_\infty E $
(cf.~[W, Sec.~4.2 ]).  On letting $z \in \Omega$ approach
$\partial_\infty E $ and appealing to the lower semi-continuity of
potentials, it follows that $\frac{1}{2\pi}\frac{\partial
g_\Omega}{\partial n}\dd s$ has energy at most $V_E$, and so by uniqueness
of the minimizing measure, we deduce that
$\dd\mu_E=\frac{1}{2\pi}\frac{\partial g_\Omega}{\partial n}\dd s$. Thus,
for any Borel subset $\gamma$ of $\partial_\infty E$,
$$\mu_E(\gamma)=\frac{1}{2\pi}\int_{\gamma}\frac{\partial
g_\Omega}{\partial n}\dd s=\frac{1}{2\pi}\int_{\gamma}|\Phi^\prime|\dd s.$$
  Alternatively, $\mu_E(\gamma)$ is given by the normalized angular
  measure of the image $\Phi(\gamma)$: \begin{equation}
\label{normalized angular}
\mu_E(\gamma)=\frac{1}{2\pi}\int_{\Phi(\gamma)}\dd \theta
 \end{equation}
(for this representation, the smoothness of $\partial_\infty E$ is not
needed).

The Green function is especially useful for estimating the modulus of a
polynomial outside $E$ when its sup norm on $E$ is known.
\newline
\newline
 {\bf Lemma \arabic{section}.\arabic{example}.} \addtocounter{example}{1}
(\dword{Bernstein-Walsh}).  {\sl If $p_n(z) $ is any
polynomial of degree $\le n$, then \begin{equation} |p_n(z)|\le \lVert
p_n\rVert_E \ \ee^{ng_\Omega (z, \infty)}, \ \ \ \ \ z\in\Omega,
\label{bernstein lemma}
 \end{equation} where $\lVert p_n\rVert_E :=
\max_{z\in E} |p_n(z)|$ and $\Omega$ is the unbounded component of
$\mathbb{C}\setminus E$.}
\newline
\newline
 {\bf Proof.}  Assume that
$p_n$ has exact degree $n$.  Then (\ref{bernstein lemma}) is equivalent
to \begin{equation} \frac{1}{n} \log\frac{1}{|p_n(z)|} + g_\Omega (z,
\infty) \ge \frac{1}{n}\log \frac{1}{\lVert p_n \rVert_E}, \ \ \ \ \
z\in\Omega.  \label{bernstein proof}
 \end{equation} Let $u(z)$ denote
the left-hand side of (\ref{bernstein proof}) and note that $u$ is
superharmonic in $\Omega$ and harmonic at $\infty$.  Moreover, from
property (c) for $g_\Omega$, we deduce that
$$\liminf_{\substack{z\rightarrow\zeta \\ z\in\Omega}} u(z) \ge
\frac{1}{n} \log  \frac{1}{\lVert p_n \rVert_E} \text{ for q.e.\  }
\zeta\in\partial\Omega.$$ Thus (\ref{bernstein proof}) follows from the
Minimum Principle for Superharmonic Functions.  $\Box$
\newline

It is useful to consider Green functions with poles at finite points of
the plane.  If $D$ is a domain with $\cp (\partial D)>0$, $\partial
D\subset\C$ compact, the \dword{Green function $g_D(z,\zeta)$ for $D$ with
pole at $\zeta\in D$} is the unique real-valued function of $z$
satisfying
\begin{enumerate}[(a')]
\item $g_D(z,\zeta)$ is harmonic in
$D\setminus\{\zeta\}$ and bounded outside any neighborhood of $\zeta$.
\item $g_D(z,\zeta) - \log\frac{1}{|z-\zeta|} = O(1)$ as
$z\rightarrow\zeta$.
\item ${\displaystyle \lim_{\substack{z\rightarrow
w \\ z\in D}}} \ g_D(z,\zeta) = 0$ for q.e.\  $w\in\partial D$.
\end{enumerate}

The relation of this Green function to the one with pole at $\infty$ is
easy to see.  Consider the mapping $w=1/(z-\zeta)$ that takes $\zeta$ to
$\infty$ and $D$ to some domain $D'$.  Then \begin{equation}
g_D(z,\zeta) = g_{D'}\left(\frac{1}{z-\zeta}, \infty\right).
\label{green finite}
 \end{equation}
\newline
\newline
 {\bf \em
Exercise.}  Verify that the function of $z$ on the right-hand side of
(\ref{green finite}) satisfies properties (a'), (b'), (c').
\newline
\newline
 {\bf \em Exercise.}  Verify that for the unit disk $D: |z|<1$,
$$g_D(z,\zeta) = \log \left|\frac{1-\overline{\zeta}z}{z-\zeta}\right|,
\ \ \ \ \ z, \zeta \in D.$$

A clever application of Green's formula shows that $g_D$ is symmetric:
$g_D(z,\zeta) = g_D(\zeta,z)$; see [Ts].
\newline
\newline
 {\bf \em
Exercise.}  Prove that if $p_n$, $n=1, 2, \ldots$, is a sequence of
monic polynomials of respective degrees $n$ that satisfy
$$\limsup_{n\rightarrow\infty}\,\lVert p_n\rVert_{[-1,1]}^{1/n}\, \le
\frac{1}{2} = \cp([-1,1]),$$ then the proportion of the number of zeros
of the $p_n$'s that lie outside any neighborhood of $[-1,1]$ tends to
zero as $n\rightarrow\infty$.  (Recall Example 2.9 and the remark
following it.) [Hint: Let $\{z_{n,k}\}_{k\in J_n}$ denote the zeros of
$p_n$ that lie at a distance $\ge\epsilon>0$ from $[-1,1]$ and consider
the functions $$\frac{1}{n}\log\frac{1}{|p_n(z)|} - \frac{1}{n}
\sum_{k\in J_n} g(z,z_{n,k}) + g(z,\infty),$$ where
$g=g_{\overline{\C}\setminus[-1,1]}$.]
\newline

Just as $\log\frac{1}{|z-t|}$ serves as the kernel for logarithmic
potential theory, so too does $g_D(z,t)$ serve as the kernel for {\em
Green potential theory}.  If $\nu$ is a finite positive measure on $D$
with compact support in $\C$, we define \begin{equation} U_D^\nu (z) :=
\int g_D(z,\zeta) \dd \nu(\zeta), \ \ \ \ \ z\in D, \label{green potential}
\end{equation} and note that $U_D^\nu \ge 0$ in $D$ and $U_D^\nu$ is
superharmonic in $D$ and harmonic in $D\setminus\mbox{supp}(\nu)$.
Furthermore, if $\nu$ has compact support in $D$, then
$$\lim_{\substack{z\rightarrow w \\ z\in D}} U_D^\nu(z) = 0 \text{ for
q.e.\  } w\in \partial D.$$ The \dword{Green energy} of a measure $\nu$ is
defined by
\begin{equation} I^D(\nu) := \int \int g_D (z,\zeta) \dd \nu(z)
\dd \nu(\zeta).  \label{green energy}
\end{equation}

For a closed subset $E\subset D$ of positive logarithmic capacity, we
consider the minimum energy problem \begin{equation} V_E^D := \inf \{
I^D(\nu) : \nu\in\M(E)\} \label{green min energy}
 \end{equation} for
which there exists a unique measure (the \dword{Green equilibrium
measure}) $\mu_E^D\in\M(E)$ such that $I^D(\mu_E^D) = V_E^D$.  Analogous
to Frostman's Theorem 1.12, there holds
\begin{align}
U_D^{\mu_E^D}(z) & = V_E^D \text{ q.e.\  on } E \label{green frostman 1}
\\ U_D^{\mu_E^D}(z) & \le  V_E^D \text{ for all } z\in D. \label{green
frostman 2}
 \end{align} The constant
\begin{equation} \cp(E,\partial D):= \frac{1}{V_E^D} \label{condenser cap}
\end{equation}
is called the \dword{capacity of the condenser $(E,\partial D)$}; see Section 5.
\newline
\newline
 {\bf Balayage}
\newline
 Let
$D\subset\overline{\C}$ be an open set with compact boundary $\partial
D$ of positive capacity and let $\mu$ be a measure with
$\mbox{supp}(\mu)\subset\overline{D}$.  The problem of \dword{balayage} (a
French word meaning ``sweeping'') consists of finding a new measure
$\mu^b$ {\em supported on} $\partial D$ such that $\mu^b(\C) = \mu(\C)$
and \begin{equation} U^{\mu^b}(z) = U^\mu(z) \text{ for q.e.\  } z\not\in
D. \label{balayage}
 \end{equation} For a bounded domain $D$, the
sweeping out of the measure $\mu$ to $\partial D$ can always be
accomplished, but if the domain $D \subset \overline{\mathbb{C}}$
contains the point at infinity, it is necessary to modify
(\ref{balayage}) so that it reads \begin{equation} U^{\mu^b}(z) =
U^\mu(z) +c \text{ for q.e.\  } z\not\in D, \label{balayage unbounded}
\end{equation} for some constant $c$.  Necessarily \begin{equation} c =
\int g_D(z,\infty) \dd \mu(z).  \label{balayage constant}
 \end{equation}
 If $D$ is connected and regular, then equality in (\ref{balayage}) and
 (\ref{balayage unbounded}) holds for all $z\not\in D$.  To ensure
 uniqueness of $\mu^b$, a condition such as boundedness of $U^{\mu^b}$
 on $\partial D$ suffices.
\newline
\newline
 {\bf \em Exercise.}
Verify that the constant $c$ in (\ref{balayage unbounded}) is given by
(\ref{balayage constant}).  [Hint: Starting with $$\int
g_D(z,\infty)\dd \mu(z) = \int[V_E-U^{\mu_{\partial D}}(z)]\dd \mu(z),$$ use
the reciprocity law together with  (\ref{balayage unbounded}) on
$\partial D$.]
\newline

Balayage measures can also be characterized by the following property:
if $h$ is any function that is continuous on $\overline{D}$ and harmonic
in $D$, then \begin{equation} \int_D h\dd\mu = \int_{\partial D}
h\dd\mu^b.  \label{balayage measure}
 \end{equation}
\newline
\newline
{\bf Example \arabic{section}.\arabic{example}.}
\addtocounter{example}{1} If $D$ is the unit disk $|z|<1$ and
$\delta_\zeta$ is the unit point mass at $\zeta\in D$, then the balayage
of $\delta_\zeta$ to $\partial D: |z|=1$ is given by
$$\dd\delta_\zeta^b(t) = \frac{1}{2\pi}P(t,\zeta)\dd\theta,\ \ \
t=\ee^{\ii\theta},\ \ \ 0\le\theta\le 2\pi,$$ where $P$ denotes the Poisson
kernel (\ref{poisson 2}).
\newline

If $\overline{D}$ is a compact subset of $\C$ with
$\cp(\overline{D})>0$, then the balayage of $\delta_\infty$ onto the
outer boundary of $D$ is the equilibrium measure $\mu_{\overline{D}}$.

The notion of balayage is intimately connected to the Green function of
a domain.  Indeed, if $D$ is bounded, $\zeta\in D$, and $\delta_\zeta^b$
denotes the balayage of $\delta_\zeta$ to $\partial_\infty D$, then
\begin{equation} g_D(z,\zeta) = \log \frac{1}{|z-\zeta|} -
U^{\delta_\zeta^b}(z), \label{green balayage}
 \end{equation} since, as
can be verified, the right-hand satisfies properties (a'), (b'), and
(c') that characterize the Green function with pole at $\zeta$.  If
$\mu$ is a finite positive measure on $D$, then \begin{equation} \mu^b =
\int\delta_\zeta^b\dd\mu(\zeta) \label{balayage mu}
 \end{equation} and
so, on integrating (\ref{green balayage}) with respect to $\dd\mu(\zeta)$,
we get \begin{equation} U_D^\mu(z) = \int g_D(z,\zeta)\dd\mu(\zeta) =
U^\mu(z) - U^{\mu^b}(z), \ \ \ z\in D. \label{balayage integral}
\end{equation}

As a consequence of (\ref{balayage integral}) and the nonnegativity of
the Green potential, we get \begin{equation} U^{\mu^b}(z) \le U^\mu(z) \
\ \ \ \ \forall z\in \C. \label{balayage inequality}
 \end{equation}

In case $D$ is an unbounded domain with $\partial D\subset \C$, then
(\ref{green balayage}), (\ref{balayage mu}), and (\ref{balayage
integral}) must be modified to include the constant $c$ of
(\ref{balayage constant}); e.g., (\ref{balayage inequality}) becomes
\begin{equation} U^{\mu^b}(z) \le U^\mu(z) +c \ \ \ \ \ \forall z\in \C.
\label{bal ineq unbounded}
 \end{equation}

\setcounter{equation}{0} \setcounter{example}{1}
\sect{Applications to Polynomial Approximation of Analytic Functions}

Let $f$  be a continuous complex-valued function on a compact set
$E\subset\C$ and let \begin{equation} \label{eq:2.1}
e_n(f;E)=e_n(f):=\min_{p\in\Po_n}\|f-p\|_E
 \end{equation} be the error
in best uniform approximation of $f$ by polynomials of degree at most
$n$. We denote by $p_n^*$ the polynomial of best approximation:
$\|f-p_n^*\|_E=e_n(f)$.

If $e_n(f)\to0$ as $n\to\infty$, the series
$$p_1^*+\sum_{n=1}^\infty(p_{n+1}^*-p_n^*)$$ converges to $f$ uniformly
on $E$,  so that the continuous function $f$ must be analytic at every
interior point of $E$. (The collection of all  functions that are
continuous on $E$ and analytic in the interior of $E$ is denoted by
$\A(E)$.) Furthermore, it follows from  the maximum principle for
analytic functions, that the above series automatically converges on
every bounded component of $\C \setminus E$, so that its sum represents
an analytic continuation of $f$  to these components (e.g., if $E$ is
the unit circle $|z|=1$, then the convergence holds in the unit disk
$|z|\leq1$).  Such a continuation, however, may be impossible.
Therefore, in order to ensure that $e_n(f)\to0$ for {\em every function}
$f$ in $\A(E)$, it is necessary to  assume that the only component of
$\C \setminus E$ is the unbounded one; that is, $\C \setminus E$  is
connected (so that $E$ does not separate the plane).

A celebrated theorem of S.N. Mergelyan (cf.~[Ga]) asserts that this
assumption is also sufficient.  Here, we prove this result in a special
case when $E$ has a connected and regular complement $\Omega :=
\overline\C \setminus E$ and $f$ is analytic in some neighborhood of
$E$.  Our aim is to determine the {\em rate of approximation}.

For any $R > 1$, let $\Gamma_R$ denote the level
curve $\{z: g_\Omega(z,\infty) = \log R\}$, see Fig.~2 (we call such a
curve a \dword{level curve with index} $R$).   The assumption that
$\Omega$ is regular ensures that for any open set $V$ containing $E$,
the level curve $\Gamma_R$ will lie in $V$ for $R$ sufficiently close to
1.
\begin{figure}[h]
\centering \input{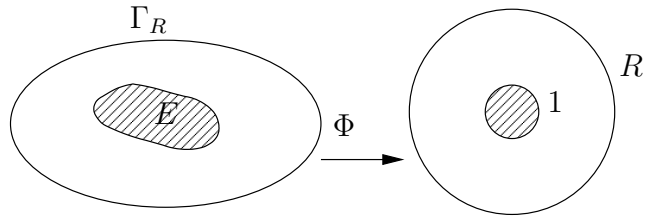} \caption{Level
curve of $g_\Omega(z,\infty)$}
 \end{figure}
\newline
 Let $F_{n+1}$ be
the $(n+1)$-st degree Fekete polynomial  for $E$ and let $P_n$ be the
polynomial of degree $\leq n$ that interpolates $f$ at the zeros of
$F_{n+1}$. We are given that $f$ is analytic in a neighborhood of $E$;
hence there exists $R > 1$ such that $f$ is analytic on and inside
$\Gamma_R$. For any such $R$, the \dword{Hermite interpolation formula}
yields \begin{equation} \label{eq:2.4} f(z)-P_n(z)=\frac{1}{2\pi
\ii}\int_{\Gamma_R}\frac{F_{n+1}(z)}{F_{n+1}(t)}\frac{f(t)\dd t}{t-z}, \;\;\;
z\mbox{ inside }\Gamma_R.
 \end{equation} (The validity of the Hermite
formula follows by first observing that the right-hand side vanishes at
the zeros of $F_{n+1}(z)$, and then by replacing $f(z)$ by its Cauchy
integral representation to deduce that the difference between $f$ and
the right-hand side is indeed a polynomial of degree at most $n$.)

Formula $(\ref{eq:2.4})$ leads to a simple estimate:
$$e_n(f)\leq\|f-P_n\|_E\leq
K\frac{\|F_{n+1}\|_E}{\min_{\Gamma_R}|F_{n+1}(t)|},$$ where $K$ is some
constant independent of $n$.  Applying parts (b), (c) of the Fundamental
Theorem 1.18, we obtain that 
\begin{equation}
\label{eq:2.5}
\limsup_{n\to\infty}e_n(f)^{1/n}\leq\frac{\cp(E)}{R\:\cp(E)}=\frac1R<1.
\end{equation} We have proved that indeed $e_n(f)\to0$ and that the
convergence is geometrically fast.  Since $R > 1$ was arbitrary (but
such that $f$ is analytic on and inside $\Gamma_R$), we have actually
proved that $(\ref{eq:2.5})$ holds with $R$ replaced by $R(f)$, where
$$R(f) := \sup \{ R : \; f \mbox{ admits analytic continuation to the
interior of } \Gamma_R\}.$$ Can we improve on this?  The answer is ---
no!  In order to show this, we appeal to the Bernstein-Walsh Lemma 3.7.

Assume now that $(\ref{eq:2.5})$ holds for some $R > R(f)$  and let
$R(f) <\rho< R$.  Then for some constant $c > 1$,
$$e_n(f)\leq\frac{c}{\rho^n}, \;\;\; n\geq1.$$ Since, from the triangle
inequality, $$\|p_{n+1}^*-p_n^*\|_E=\|p_{n+1}^*-f+f-p_n^*\|_E\leq
e_{n+1}(f)+e_n(f)\leq2c\rho^{-n},$$ we obtain from the Bernstein-Walsh
Lemma that for any $r > 1$,
$$\|p_{n+1}^*-p_n^*\|_{\Gamma_R}\leq2c\left(\frac{r}{\rho}\right)^n,
\;\;\; n\geq1.$$ If we choose $R(f ) < r < \rho$, we obtain that the
series $p_1^*+\sum_{n=1}^\infty(p_{n+1}^*-p_n^*)$ converges uniformly
inside $\Gamma_r$. Hence it gives an analytic continuation of $f$ to the
interior of $\Gamma_r$, which contradicts the definition of $R(f)$.

Let us summarize what we have proved.
\newline
\newline
 {\bf Theorem \arabic{section}.\arabic{example}} \addtocounter{example}{1}
(\dword{Walsh} [W, Ch.$\:$VII]).  {\sl Let $E$ be a compact set with connected
and regular complement.  Then for any $f\in\A(E)$,
$$\limsup_{n\to\infty}e_n(f)^{1/n}=\frac{1}{R(f)}.$$}
\newline
\newline
{\em Remark.} $R(f)$ is the first value of $R$ for which the level curve
$\Gamma_R$ contains a singularity of $f$. It may well be possible that
$f$ is analytic at some other points of $\Gamma_{R(f)}$, but the
geometric rate of best polynomial approximation ``does not feel this''
--- whether every point of $\Gamma_{R(f)}$ is a singularity or merely
one point is a singularity, the rate of approximation remains the same
as if  $f$  was analytic only inside of $\Gamma_{R(f)}!$ To take
advantage of any extra analyticity, different approximation tools are
needed; e.g., rational functions.
\newline
\newline
 {\bf Example \arabic{section}.\arabic{example}.} \addtocounter{example}{1} Let
$E=[-1,-\alpha]\cup[\alpha,1]$, $0<\alpha<1$, and let $f=0$ on
$[-1,-\alpha]$ and $f=1$ on $[\alpha,1]$. Some level curves $\Gamma_R$
of $g_{\C\setminus E}$ are depicted on Fig.~3.
For $R$ small, $\Gamma_R$ consists of two pieces, while for $R$ large,
$\Gamma_R$ is a single curve.  There is a ``critical value''
$R_0=g_{\C\setminus E}(0,\infty)$ for which $\Gamma_{R_0}$ represents a
self-intersecting lemniscate-like curve (the bold curve in Fig.~3).
Clearly, $f$ can be extended as an analytic function
to the interior of $\Gamma_{R_0}$ (define $f=0$ inside the left lobe and
$f=1$ inside the right lobe).  For $R>R_0$, the interior of $\Gamma_R$
is a (connected) domain; hence there is no function analytic inside of
$\Gamma_R$ that is equal to 0 on $[-1,-\alpha]$ and to 1 on
$[\alpha,1]$. Therefore $$R(f)=R_0=\exp\left\{g_{\C\setminus
E}(0,\infty)\right\},$$ and by Theorem 4.1:
$$\limsup_{n\to\infty}e_n(f)^{1/n}=\exp\left\{-g_{\C\setminus
E}(0,\infty)\right\}.$$
\begin{figure}[h]
\centering \input{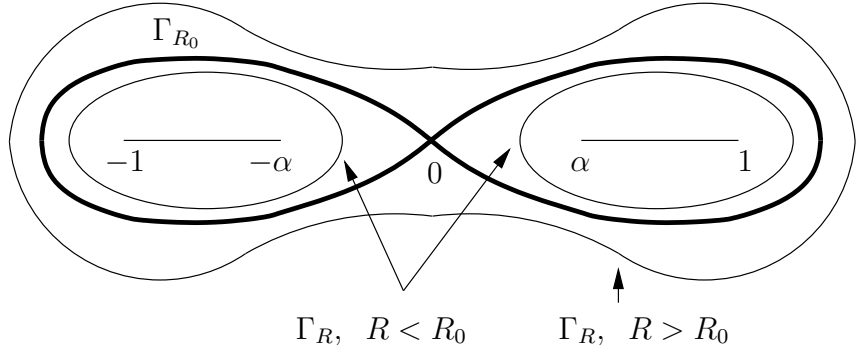}
\caption{Level curves of $g_{\C\setminus E}$ for Example 4.2}
\end{figure}

\setcounter{equation}{0} \setcounter{example}{1}

\sect{Rational Approximation}

For a rational function $R (z) = P_1 (z) / P_2 (z)$, where $P_1$ and
$P_2$ are monic polynomials of degree $n$, one can write
$$-\frac1n\log|R(z)|=U^{\nu_1}(z)-U^{\nu_2}(z),$$ where $\nu_1, \nu_2$
are the normalized zero counting measures for $P_1$, $P_2$,
respectively.  The right-hand side represents the logarithmic potential
of the {\em signed measure} $\mu=\nu_1-\nu_2$:
$$U^{\nu_1}(z)-U^{\nu_2}(z)=U^\mu(z)=\int\log\frac{1}{|z-t|}\dd \mu(t).$$
The theory of such potentials can be developed along the same lines as
in the earlier sections.  We present below only the very basic notions
of this theory that are needed to formulate the approximation results.
A more in-depth treatment can be found in the works of Bagby [B],
Gonchar [Gon], as well as [ST].

The analogy with electrostatics problems suggests considering the
following energy problem.  Let $E_1, E_2\subset\C$  be two closed sets
that are a  positive distance apart.  The pair $(E_1,E_2)$ is called a
\dword{condenser} and the sets $E_1$,  $E_2$ are called the \dword{plates.}
Let $\mu_1$ and $\mu_2$  be positive unit measures supported on $E_1$
and $E_2$, respectively.  Consider the energy integral of the signed
measure $\mu=\mu_1-\mu_2$:
$$I(\mu)=\int\int\log\frac{1}{|z-t|}\dd \mu(z)\dd \mu(t).$$ Since   $\mu(\C) =
0$, the integral is well-defined, even if one of the sets is unbounded.
While not obvious, it turns out that such $I(\mu)$ is always {\em
positive.} We assume that $E_1$  and $E_2$  have positive logarithmic
capacity.  Then the minimal energy (over all signed measures of the
above form) $$V(E_1,E_2):=\inf_\mu I(\mu)$$ is finite and positive.  We
then define the \dword{condenser capacity} $\cp(E_1,E_2)$ by
$$\cp(E_1,E_2):=1/V(E_1,E_2).$$ One can show, as with the Frostman
theorem, that there exists a unique signed measure \linebreak
$\mu^*=\mu_1^*-\mu_2^*$ (the {\em equilibrium measure} for the
condenser) for which $I(\mu^*)=V(E_1,E_2)$. Furthermore, the
corresponding potential (called the \dword{condenser potential}) is
constant on each plate:
\begin{equation} \label{eq:6.1} U^{\mu^*}=c_1
\mbox{ on } E_1, \;\;\; U^{\mu^*}=-c_2 \mbox{ on } E_2
\end{equation}
(we assume throughout that  $E_1$, $E_2$  are regular --- otherwise the
above equalities hold only quasi-everywhere).  On integrating against
$\mu^*$, we deduce from $(\ref{eq:6.1})$ that \begin{equation}
\label{eq:6.2} c_1+c_2=V(E_1,E_2)=1/\cp(E_1,E_2).
  \end{equation} We
mention that (similar to the case of the conductor potential) the
relations of type $(\ref{eq:6.1})$ {\em characterize} $\mu^*$. Moreover,
one can deduce from $(\ref{eq:6.1})$ that the measure $\mu_i^*$ is
supported on the boundary (not necessarily the outer one) of $E_i$, $i =
1,2$.  Therefore, on replacing each $E_i$ by its boundary, we do not
change the condenser capacity or the condenser potential.
\newline
\newline
 {\bf Example \arabic{section}.\arabic{example}.}
\addtocounter{example}{1}  Let $E_1$, $E_2$ be, respectively,  the
circles  $|z| = r_1$ $|z| = r_2$, $r_1<r_2$. These sets are invariant
under rotations.  Being unique, the measure $\mu^*$ is therefore also
invariant under rotations and we obtain  that $$\dd\mu_1^*=\frac{1}{2\pi
r_1}\dd s, \;\;\; \dd \mu_2^*=\frac{1}{2\pi r_2}\dd s,$$ where $\dd s$ denotes
the arclength over the respective circles $E_1$, $E_2$. Applying the result
of Example 1.10, we find that
$$U^{\mu^*}(z)=\left\{
\begin{array}{ll} 0,             & |z|>r_2 \\ \log(r_2/|z|), &
r_1\leq|z|\leq r_2 \\ \log(r_2/r_1), & |z|<r_1.
  \end{array} \right.$$
Therefore (recall ($\ref{eq:6.2}$)) \begin{equation} \label{eq:6.3}
\cp(E_1,E_2)=1/\log\frac{r_2}{r_1}.
  \end{equation}

Assume now that each plate of a condenser is a single Jordan arc or
curve (without self-intersections), and let $G$ be the doubly-connected
domain that is bounded by $E_1$ and $E_2$, see Fig.~4.  We call such a
$G$ a \dword{ring domain}.
\begin{figure}[h]
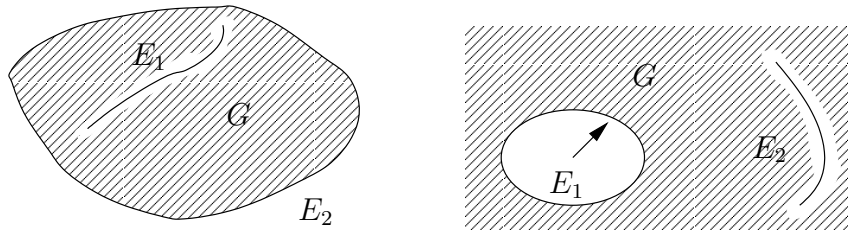

\centering \input{ptt_ar3.tex} \hspace{1cm} \input{ptt_ar4.tex}
\caption{Ring domains}
\end{figure}
For ring domains one can give an
alternative definition of condenser capacity.  Let $$u(z) := \int\log(z
- t)\dd \mu^*(t) + c_1.$$ The complex function $u$ is locally analytic but
not single-valued in $G$ (notice that there is no modulus sign in the
integral).  Moreover, if we fix $t$ and let $z$ move along a simple
closed counterclockwise oriented curve in $G$ that encircles $E_1$, say,
then the imaginary part of $\log(z - t)$ increases by  $2\pi$, for $t\in
E_1$, while for $t\in E_2$ it returns to the original value.  Since
$\mu_1^*$ and $\mu_2^*$ are unit measures, it follows that the function
$\phi: z \to  w= \exp (u(z))$ is analytic and single-valued.  Moreover,
it can be shown to be one-to-one in $G$. By its definition, $\phi$
satisfies $$\log|\phi|=-U^{\mu^*}+c_1=0 \mbox{ on } E_1; \;\;\;
\log|\phi|=-U^{\mu^*}+c_1=c_1+c_2 \mbox{ on } E_2.$$ Therefore $\phi$
maps $G$ conformally onto the annulus $1<|w|<\ee^{c_1+c_2}$.

It is known from the theory of conformal mapping that, for a ring domain
$G$, there exists unique $R > 1$, called the \dword{modulus of} $G$  (we
denote it by mod$(G)$), such that $G$ can be mapped conformally onto the
annulus $1 < |w| < R$. We have thus shown that \begin{equation}
\label{eq:6.4} \cp (E_1,E_2) = 1/\log (\mbox{mod}(G)).
  \end{equation}
We remark that if  $G_1\supset G_2$ are two ring domains, then
mod$(G_1)\geq$ mod$(G_2)$.
\newline
\newline
 {\bf Example
\arabic{section}.\arabic{example}.} \addtocounter{example}{1} Let $E_1$,
$E_2$ be as above, and assume that $E_2$  is the $R$-th level curve for
$E_1$. That is, $|\Phi(z)|=R$ for $z \in E_2$, where $\Phi$ maps
conformally the unbounded component of $\C \setminus E_1$  onto $|w| >
1$. In particular, $\Phi$ maps the corresponding ring domain $G$ onto
the annulus $1 <|w|< R$, and we conclude that mod$(G)=R$ (so that
$\cp(E_1,E_2)=1/\log R$).  Applying this to the configuration of
Example 5.1, we see that $\Phi(z) = z/r_1$, so that
$R = r_2/r_1$,  and  we obtain again $(\ref{eq:6.3})$.
\newline

We now turn to rational approximation.  Let $E\subset\C$ be compact.  We
denote by ${\cal R}_n$ the collection of all rational functions of the
form $R = P / Q$, where $P$, $Q$ are polynomials of degree at most $n$,
and $Q$ has no zeros in $E$. For  $f\in{\cal A}(E)$, let
$$r_n(f;E)=r_n(f):=\inf_{r\in{\cal R}_n}\|f-r\|_E$$ be the error in best
approximation of $f$  by rational functions from ${\cal R}_n$.  Clearly,
since polynomials are rational functions, we have
(cf.~$(\ref{eq:2.1})$) $r_n(f)\leq e_n(f)$. A basic theorem regarding the
rate of rational approximation was proved by Walsh [W, Ch.IX].
Following is a special case of this theorem.
\newline
\newline
\newline
{\bf Theorem  \arabic{section}.\arabic{example}.}
\addtocounter{example}{1} (Walsh) {\sl Let $E$ be a single Jordan arc or
curve and let $f$ be analytic on a simply connected domain $D \supset
E$. Then}
\begin{equation} \label{eq:6.5} \limsup_{n\to\infty}
r_n(f)^{1/n}\leq \exp(-1/\cp(E,\partial D)).
\end{equation}
\newline

The proof  of $(\ref{eq:6.5})$ follows the same ideas as the proof of
inequality $(\ref{eq:2.5})$.  Let $\Gamma$ be a contour in $D\setminus
E$ that is arbitrarily close to $\partial D$. Let
$\mu^*=\mu_1^*-\mu_2^*$ be the equilibrium measure for the condenser
$(E,\Gamma)$. For any $n$, let $\alpha_1^{(n)},\ldots,\alpha_n^{(n)}$ be
equally spaced on $E$ (with respect to $\mu_1^*$) and let
$\beta_1^{(n)},\ldots,\beta_n^{(n)}$ be equally spaced on $\Gamma$ (with
respect to $\mu_2^*$).  Then one can show that the rational functions
$r_n(z)$ with zeros at the $\alpha_i^{(n)}$'s and poles at the
$\beta_i^{(n)}$'s satisfy
\begin{equation} \label{eq:6.6}
\left(\frac{\displaystyle \max_E|r_n|}{\displaystyle
\min_\Gamma|r_n|}\right)^{1/n}\to \ee^{-1/\cp(E, \Gamma)}.
\end{equation}
Let $R_n=p_{n-1}/q_n$ be the rational function with poles at the
$\beta_i^{(n)}$'s that interpolates $f$ at the points
$\alpha_i^{(n)}$'s.  Then the Hermite formula (cf.~$(\ref{eq:2.4})$)
takes the following form: $$f(z)-R_n(z)=\frac{1}{2\pi
\ii}\int_\Gamma\frac{r_n(z)}{r_n(t)}\frac{f(t)}{t-z}\dd t, \;\;\; z \mbox{
inside }\Gamma,$$ and it follows from $(\ref{eq:6.6})$ that
$$\limsup_{n\to\infty}
r_n(f)^{1/n}\leq\limsup_{n\to\infty}\|f-R_n\|_E^{1/n}\leq \ee^{-1/\cp(E,
\Gamma)}.$$ Letting $\Gamma$ approach $\partial D$, we get the result.
\newline
\newline
 {\bf Remarks.}
\newline
 (a) Unlike in the polynomial
approximation, no rate of convergence of $r_n(f)$ to 0 can ensure that a
function  $f\in C(E)$ is analytic somewhere beyond $E$.
\newline
 (b) One
can construct a function for which equality holds in $(\ref{eq:6.5})$,
so that this bound is sharp.  Such a function necessarily has a
singularity at every point of $\partial D$; otherwise $f$ would be
analytic in a larger domain, so that the corresponding condenser
capacity will become smaller.  In view of Theorem 5.3, this would
violate the assumed equality in $(\ref{eq:6.5})$.
\newline

Although sharp, the bound $(\ref{eq:6.5})$ is unsatisfactory, in the
following sense.  Assume, for example, that $E$ is connected and has a
connected complement, and let $\Gamma_R$, $R>1$, be a level curve for
$E$. Let $f$ be a function that is analytic in the domain $D$ bounded by
$\Gamma_R$ and such that the equality holds in $(\ref{eq:6.5})$.
According to Example 5.2, we then obtain that
\begin{equation} \label{eq:6.7}
\limsup_{n\to\infty}r_n(f)^{1/n}=\frac1R.
 \end{equation} By Remark (b)
above, such $f$ must have singularities on $\Gamma_R$. Hence (recall
Remark following
Theorem 4.1) the relation
$(\ref{eq:6.7})$ holds with $r_n(f)$ replaced by $e_n(f)$. But the
family ${\cal R}_n$ contains $\Po_n$ and it is much more rich than
$\Po_n$ --- it depends on $2n+1$ parameters while $\Po_n$ depends only
on $n+1$ parameters.  One would expect, therefore, that at least for a
subsequence of $n$'s, $r_n(f)$ behaves asymptotically like $e_{2n}(f)$.
This was a motivation for the following conjecture.
\newline
\newline
{\bf Conjecture.} (A.A. Gonchar) {\sl Let $E$ be a compact set and $f$
be analytic in an open set $D$ containing $E$.  Then} \begin{equation}
\label{eq:6.8} \liminf_{n\to\infty}r_n(f;E)^{1/n}\leq
\exp(-2/\cp(E,\partial D)).
  \end{equation} This conjecture was proved
by O. Parfenov [Pa] for the case when $E$ is a continuum with connected
complement and in the general case by V. Prokhorov [P]; they used a very
different method --- the so-called ``AAK Theory'' (cf.~[Y]).  However
this method is not constructive, and it remains a challenging problem to
find such a method.  Yet, potential theory can be used to obtain bounds
like $(\ref{eq:6.8})$ in the stronger form
$$\lim_{n\to\infty}r_n(f;E)^{1/n}= \exp(-2/\cp(E,\partial D))$$ for
some important {\em subclasses} of analytic functions, such as Markov
functions (cf.\ [Gon]) and functions with a finite number of algebraic
branch-points (cf.\ [St]).

\setcounter{equation}{0} \setcounter{example}{1}
\sect{Logarithmic Potentials with External Fields}

Let $E$ be a closed (not necessarily compact) subset of $\C$ and let
$w(z)$  be a nonnegative weight on $E$. We define a new ``distance
function'' on $E$, replacing  $|z-t|$  by  $|z-t|w(z)w(t)$. This gives
rise to weighted versions of logarithmic capacity, transfinite diameter
and Chebyshev constant.
\newline
\newline
 {\bf Weighted capacity:} $\cp(w, E)$.

As before, let $\M(E)$ denote  the collection of all unit measures
supported on $E$. We set $$Q:=\log\frac1w$$ and call it the
\dword{external field}.  Consider the modified energy integral for
$\mu\in\M(E)$:
\begin{eqnarray} \label{eq:4.1} I_w(\mu) &:=&
\int\int\log\frac{1}{|z-t|w(z)w(t)}\dd \mu(z)\dd \mu(t)  \\ {}       &=&
\int\int\log\frac{1}{|z-t|}\dd\mu(z)\dd \mu(t)+2\int Q(z)\dd \mu(z) \nonumber
\end{eqnarray}
and let $$V_w:=\inf_{\mu\in\M(E)}I_w(\mu).$$ The
\dword{weighted capacity} is defined  by $$\cp(w,E):=\ee^{-V_w}.$$ In the sequel,
we assume that $w$ satisfies the following conditions:
\begin{itemize}
  \item[(i)]   $w > 0$  on a subset of positive logarithmic capacity;
  \item[(ii)]  $w$ is continuous (or, more generally, upper
  semi-continuous);
  \item[(iii)] if  $E$ is unbounded, then
  $|z|w(z)\to0$ as $|z|\to\infty$, $z\in E$.
 \end{itemize}

Under these restrictions on $w$, there exists a unique measure
$\mu_w\in\M(E)$, called the \dword{weighted equilibrium measure}, such
that $$I(\mu_w)=V_w.$$ The above integral $(\ref{eq:4.1})$ can be
interpreted as the total energy of the unit charge $\mu$, in the
presence of the external field $Q$ (in this electrostatics
interpretation, the field is actually $2Q$).  Since this field has a
strong repelling effect near points where $w = 0$ (i.e.,  $Q=\infty$),
assumption (iii) physically means that, for the equilibrium
distribution, no charge occurs near $\infty$. In other words, the
support $\mbox{supp}(\mu_w)$ of $\mu_w$ is necessarily {\em compact}.
However, unlike the unweighted case, the support need not lie entirely
on $\partial_\infty E$ and, in fact, it can be quite an arbitrary closed
subset of $E$.  Determining this set is one of the most important
aspects of weighted potential theory.
\newline
\newline
 {\bf Weighted
transfinite diameter:}  $\tau(w, E)$.

Let $$\delta_n(w):=\max_{z_1,\ldots,z_n\in E}\left(\prod_{1\leq i<j\leq
n}|z_i-z_j|w(z_i)w(z_j)\right)^{2/n(n-1)}.$$ Points
$z_1^{(n)},\ldots,z_n^{(n)}$ at which the maximum is attained are called
\dword{weighted Fekete points}.  The corresponding
\dword{Fekete polynomial}
is the monic polynomial with all its zeros at these points.

As in the unweighted case, the sequence $\delta_n(w)$ is decreasing, so
one can define $$\tau(w,E):=\lim_{n\to\infty}\delta_n(w),$$ which we
call the \dword{weighted transfinite diameter} of $E$.
\newline
\newline
{\bf Weighted Chebyshev constant:} $\cheb(w, E)$.

Let $$t_n(w):=\min_{p\in\Po_{n-1}}\|w^n(z)(z^n-p(z))\|_E.$$ Then the
\dword{weighted Chebyshev constant} is defined by
$$\cheb(w,E):=\lim_{n\to\infty}t_n(w)^{1/n}.$$

The following theorem generalizes the fundamental results stated in
Theorem 1.18.
\newline
\newline
 {\bf Theorem
\arabic{section}.\arabic{example}} \addtocounter{example}{1}
(\dword{Generalized Fundamental Theorem}).  {\sl Let $E$ be a closed set of
positive capacity.  Assume that $w$ satisfies the conditions
\textnormal{(i)--(iii)} and let $Q=\log(1/w)$. Then
$$\cp(w,E)=\tau(w,E)=\cheb(w,E)\exp\left(-\int Q\dd \mu_w\right).$$
Moreover, weighted Fekete points have asymptotic distribution $\mu_w$ as
$n\to\infty$, and weighted Fekete polynomials are asymptotically optimal
for the weighted Chebyshev problem.}
\newline

How can one find $\mu_w?$

In most applications, the weight $w$ is continuous and the set $E$ is
regular.  Recall that the latter means that the classical (unweighted)
equilibrium potential for $E$ is equal to $V_E$  {\em everywhere} on
$E$,  not just quasi-everywhere.  Under these assumptions, the
equilibrium measure $\mu=\mu_w$ is characterized by the conditions that
$\mu\in\M(E)$, $I(\mu)<\infty$ and, for some constant $c_w$, the
following {\em variational conditions} hold: \begin{equation}
\label{eq:4.2} \left\{ \begin{array}{ll} U^\mu+Q=c_w & \mbox{on } S(\mu)
= \mbox{supp}(\mu) \\ U^\mu+Q\geq c_w & \mbox{on } E.
 \end{array}
\right.
  \end{equation} On integrating (against $\mu=\mu_w$) the first
condition, we obtain that the constant  is given by $$c_w=I(\mu_w)+\int
Q\dd \mu_w=V_w-\int Q\dd \mu_w.$$ When trying to find $\mu_w$, an essential
step (and a nontrivial problem in its own right!) is to determine the
support $S(\mu_w) : = \mbox{supp}(\mu_w)$. There are several methods by
which $S(\mu_w)$ can be numerically approximated, but they are
complicated from the computational point of view.  Therefore, knowing
properties of the support can be useful and we list some of them.
\newline
\newline
 {\bf Properties of the support $S(\mu_w)$}
\newline
 (a) The sup norm of weighted polynomials ``lives'' on
$S(\mu_w)$. That is, for any $n$ and for any polynomial $P_n$ of degree
at most $n$, there holds $$\|w^nP_n\|_E=\|w^nP_n\|_{S(\mu_w)}.$$ (b) Let
$K$ be a compact subset of $E$ of positive capacity, and define $$F(K)
:= \log\cp(K) - \int_KQ\dd \mu_K,$$ where $\mu_K$ is the classical
(unweighted) equilibrium measure for $K$. This so-called
\dword{F-functional} of Mhaskar and Saff is often a helpful tool in finding
$S(\mu_w)$. Since $\cp(K)$ and $\mu_K$ remain the same if we replace $K$
by  $\partial_\infty K$, we obtain that $F(K)=F(\partial_\infty K)$. It
turns out that the outer boundary of $S(\mu_w)$ maximizes the
F-functional: $$\max_K F(K)=F(\partial_\infty S(\mu_K)).$$ This result
is especially useful when $E$ is a real interval and $Q$ is {\em
convex}.  It is then easy to derive from $(\ref{eq:4.2})$ that
$S(\mu_w)$ is an {\em interval}.  Thus, to find the support, one merely
needs to maximize $F(K)$ only over intervals $K\subset E$, which amounts
to a standard calculus problem  for the determination of the endpoints
of $S(\mu_w)$.
\newline
\newline
 (c) $S(\mu_w)$ is the set of weighted
polynomial \dword{peaking points}; that is, if $w$ is continuous and $E$
is of positive capacity at each of its points, then $z$ belongs to
$S(\mu_w)$ iff for every disk $D_r(z)$ there is a  weighted polynomial
$w^nP_n$ that attains its maximum modulus only in $D_r(z)$ (cf.\ [ST,
Sec.  IV.1]).
\newline
\newline
 {\bf Example
\arabic{section}.\arabic{example}.  Incomplete polynomials}
\addtocounter{example}{1}

For the study of incomplete polynomials of type $\theta$ on the interval
$E=[0,1]$; that is, polynomials of the form $p(x) = \sum_{k=s}^n a_k
x^k$ where $s/n \ge \theta$, the appropriate external field is
$Q(x)=\log(1/w(x))=-\frac{\theta}{1-\theta}\log x$ which is convex.
Maximizing the F-functional one gets $S(\mu_w)=[\theta^2,1]$. (For
details, see [ST, Sec.  IV.1].)
\newline
\newline
 {\bf Example
\arabic{section}.\arabic{example}.  Freud
Weights}\addtocounter{example}{1}

Here $E = \R$ and $w(x) = \exp ( -|x|^\alpha)$. Hence $Q(x) =
|x|^\alpha$ is convex provided that $\alpha>1$, and we obtain
$S_w=[-a_\alpha,a_\alpha]$, where $a_\alpha$ can be given explicitly in
terms of the Gamma function.  (Actually, this result also holds for all
$\alpha>0$; see [ST, Sec.  IV.1].) For example, when $\alpha=2$, we get
$S_w=[-1,1]$.
\newline

The {\em Generalized Weierstrass Approximation Problem} mentioned in
problem (v) of the introduction states the following: For $E\subset\R$
closed, $w:E\to[0,\infty)$, characterize those functions $f$ continuous
on $E$ that are uniform limits on $E$ of some sequence of weighted
polynomials $(w^nP_n)$, $\deg P_n\leq n$.

To attack this problem, we begin with a crucial observation.  Let $E$ be
a closed subset of $\R$ whose complement is regular and $w(x)$ be
continuous on $E$. Then we have the following weighted analogue of the
Bernstein-Walsh lemma:
$$|w^n(x)P_n(x)|\leq\|w^nP_n\|_{S(\mu_w)}\exp(-n(U^{\mu_w}(x)+Q(x)-c_w)),
\;\;\; x\in E\setminus S(\mu_w).$$ With the aid of $(\ref{eq:4.2})$ and
a variant of the Stone-Weierstrass theorem (cf.\ [ST]), one can show
that if a sequence $(w^n(x)P_n(x))$, $\deg P_n\leq n$, converges
uniformly on $E$, then it tends to $0$ for every $x\in E\setminus
S(\mu_w)$.

Thus, if some $f \in C(E)$ is a uniform limit on $E$ of such a sequence
as $n \to \infty$, it must vanish on $E \setminus S(\mu_w)$. The
converse is not true, in general, but it is true in many important
cases, such as for incomplete polynomials where the weight $w(x) =
x^{\theta/(1-\theta)}$ on $[0,1]$ and for Freud weights $w(x) =
\exp(-|x|^\alpha),\ \alpha>1$, on $\R$.  The latter fact provided an
essential ingredient in resolving problem (iv) of the introduction (see
[LuSa] and [LMS]).

For the case when $E$ is a real interval and $Q=\log(1/w)$ is convex on
$E$, this author conjectured and   Totik [To] has proved that, more
generally, any $f\in C(E)$ that vanishes on $E\setminus S(\mu_w)$ is the
uniform limit on $E$ of some sequence of weighted polynomials $(w^n
P_n)$, $\deg P_n\le n$.

\def\thesection{}
\sect{\hskip-.35truecm References}
\begin{itemize}

  \item [{[B]}]    T. Bagby, The modulus of a plane condenser, {\em J. Math.
Mech.}, {\bf 17} (1976) 315--329.

  \item [{[F]}]    G. Freud, On the coefficients in the recursion formulae of
orthogonal polynomials, {\em Proc.~Roy.\ Irish Acad.~Sect.~A(1)}, {\bf 76}
(1976) 1--6.

  \item [{[Ga]}]   D. Gaier, {\em Lectures on Complex Approximation},
Birkh\"{a}user, Boston Inc., Boston, MA, 1987.

   \item [{[Gon]}]  A.A. Gonchar, The rate of rational approximation of 
analytic functions. (Russian) Modern problems of mathematics. Differential 
equations, mathematical analysis and their applications.  {\em Trudy Mat. 
Inst. Steklov}, {\bf 166} (1984), 52--60.

  \item [{[La]}]  N.S. Landkof, {\em Foundations of Modern Potential Theory},
Springer-Verlag, Heidelberg, 1972.

  \item [{[Lo]}]    G.G. Lorentz, Approximation by incomplete polynomials
(problems and results). In E.B. Saff and R.S. Varga, editors, {\em Pad\'{e}
and Rational Approximations: Theory and Applications}, 289--302, Academic
Press, New York, 1977.

  \item [{[LuSa]}] D.S. Lubinsky and E.B. Saff, Uniform and mean approximation
by certain weighted polynomials, with applications, {\em Constr.~Approx.},
{\bf 4} (1988) 21--64.

  \item [{[LMS]}]  D.S. Lubinsky, H.N. Mhaskar, and E.B. Saff, Freud's
conjecture for exponential weights, {\em Bull.~Amer.~Math.~Soc.}, {\bf 15}
(1986) 217--221.

  \item [{[Pa]}]   O.G. Parf\"enov, Estimates of singular numbers of the 
  Carleson embedding operator,  
{\em Mat.~Sb.\ (N.S.)}, {\bf 131(173)} (1986) 501--518;  English transl.\ in
{\em Math.~USSR-Sb.}, {\bf 59}  (1988)  497--514.

  \item [{[P]}]    V.A. Prokhorov, Rational approximation of analytic
functions, {\em Mat.~Sb}, {\bf 184} (1993) 3--32. English transl.\ {\em Russian
Acad.~Sci.~Sb.~Math.} {\bf 78} (1994) 139--164.

  \item [{[R1]}]    T. Ransford, {\em Potential Theory in the Complex Plane},
Cambridge University Press, Cambridge, 1995.

  \item [{[R2]}]    T. Ransford, Computation of logarithmic capacity, {\em Comp.\ Methods Function Theory}, (to appear).

  \item [{[ST]}] E.B. Saff and V. Totik, {\em Logarithmic Potentials with
External Fields}, Springer-Verlag, New York, 1997.

  \item [{[St]}]   H. Stahl, General convergence results for rational
approximation, in: {\em Approximation Theory VI}, volume 2,  C.K. Chui, L.L.
Schumaker and J.D. Ward (eds.), Academic Press, Boston, (1989) 605--634.

  \item [{[Sz]}]   G. Szeg\H{o}, {\em Orthogonal Polynomials}, volume 23 of
{\em Colloquium Publications}, Amer.~Math.~Soc., Providence, R.I., 1975.

  \item [{[To]}]   V. Totik, Weighted polynomial approximation for convex
external fields, {\em Constr.~Approx.}, {\bf 16} (2000) 261--281.

  \item [{[Ts]}]   M. Tsuji, {\em Potential Theory in Modern Function Theory},
Maruzen, Tokyo, 1959.

  \item [{[W]}]    J.L. Walsh, {\em Interpolation and Approximation by Rational
Functions in the Complex Plane}, volume 20 of {\em Colloquium Publications,}
Amer.~Math.~Soc., Providence, R.I., 1960.

  \item [{[Y]}]    N. Young, {\em An Introduction to Hilbert Space}, Cambridge
University Press, Cambridge, 1988.

\end{itemize}



{

\bigskip
\hskip1.4 em\vbox{\noindent Edward B. Saff
 \\ Center for Constructive Approximation \\ Department of Mathematics
 \\ Vanderbilt University \\ Nashville TN 37240 USA \\
 {\tt edward.b.saff@vanderbilt.edu}\\
{\tt http://www.math.vanderbilt.edu/\~{}esaff}}

}

\endddoc